\documentclass[10pt]{amsart}
\usepackage{graphicx} 
\usepackage{multicol} 
\usepackage{amsmath,amsthm,amssymb,amsfonts, eucal, amscd, mathbbol, mathrsfs, mathabx}
\usepackage[shortlabels]{enumitem}
\usepackage{cite}
\usepackage{setspace}
\usepackage[all,cmtip]{xy}{
\usepackage{graphicx}
\usepackage{subfig}
\usepackage{fancyhdr}
\usepackage{latexsym}
\usepackage{fncylab}
\usepackage{epic}
\usepackage{ifthen}
\usepackage{xcolor}
\usepackage{amscd}
\usepackage{indentfirst}
\usepackage{cite}
\usepackage{subfig}
\usepackage{rotating}
\usepackage{scalefnt}
\usepackage{enumitem}
\setlist{nolistsep}

\parskip = 0.2in
\parindent = 0.0in
\topmargin = 0.0in

\oddsidemargin = 0.0in
\evensidemargin = 0.0in
\textwidth = 6.0in

\newtheorem{thm}{thm}[section]
\newtheorem*{thm*}{Main Theorem}
\newtheorem*{thm**}{thm}

\newtheorem{cor}[thm]{Corollary}

\theoremstyle{defn}
\newtheorem{defn}[thm]{Definition}
\theoremstyle{defn}

\theoremstyle{defn}

\theoremstyle{defn}

\theoremstyle{defn}

\theoremstyle{defn}

\theoremstyle{defn}

\numberwithin{thm}{subsection}

\newcommand{\R}{\ensuremath{\mathbb{R}}}

\newcommand{\Z}{{\Bbb Z}} 

\def\p{\partial}
\def\i{\infty}
 
\def\supp{\it{supp}}

\def\cal{\mathcal}

\def\g{\gamma} 
\def\G{\Gamma}
\def\d{\delta}
 
\def\e{\epsilon}

\def\L{\Lambda}

\def\o{\omega} 
\def\O{\Omega}

\def\H{\mathcal{H}}

\def\Gr{\mathrm{Gr}}

\begin{document}
 
		\author[J. Harrison \& H. Pugh]{J. Harrison \\Department of Mathematics \\University of California, Berkeley  \\ H. Pugh\\ Mathematics Department \\ Stony Brook University} 
		\title[Plateau's Problem]{Plateau's Problem: What's Next}

\begin{abstract}Plateau's problem is not a single conjecture or theorem, but rather an abstract framework, encompassing a number of different problems in several related areas of mathematics. In its most general form, Plateau's problem is to find an element of a given collection \( \cal{C} \) of ``surfaces'' specified by some boundary constraint, which minimizes, or is a critical point of, a given ``area'' function \( F:\cal{C}\to \R \). In addition, one should also show that any such element satisfies some sort of regularity, that it be a sufficiently smooth manifold away from a well-behaved singular set. The choices apparent in making this question precise lead to a great many different versions of the problem. Plateau's problem has generated a large number of papers, inspired new fields of mathematics, and given rise to techniques which have proved useful in applications further afield. In this review we discuss a few highlights from the past hundred years, with special attention to papers of Federer, Fleming, Reifenberg and Almgren from the 1960's, and works by several groups, including ourselves, who have made significant progress on different aspects of the problem in recent years. A number of open problems are presented.
	\end{abstract}
	
	\maketitle


\section{Introduction}

Plateau's problem has intrigued mathematicians and scientists alike for over two hundred years. It remains one of the most accessible problems in mathematics, yet retains a subtle difficulty in its formulation. Many different versions of Plateau's problem have been solved, but even today there are still important questions left unanswered, and deep mysteries about the problem still remain.

\begin{figure}
	\centering
	\includegraphics[scale=0.5]{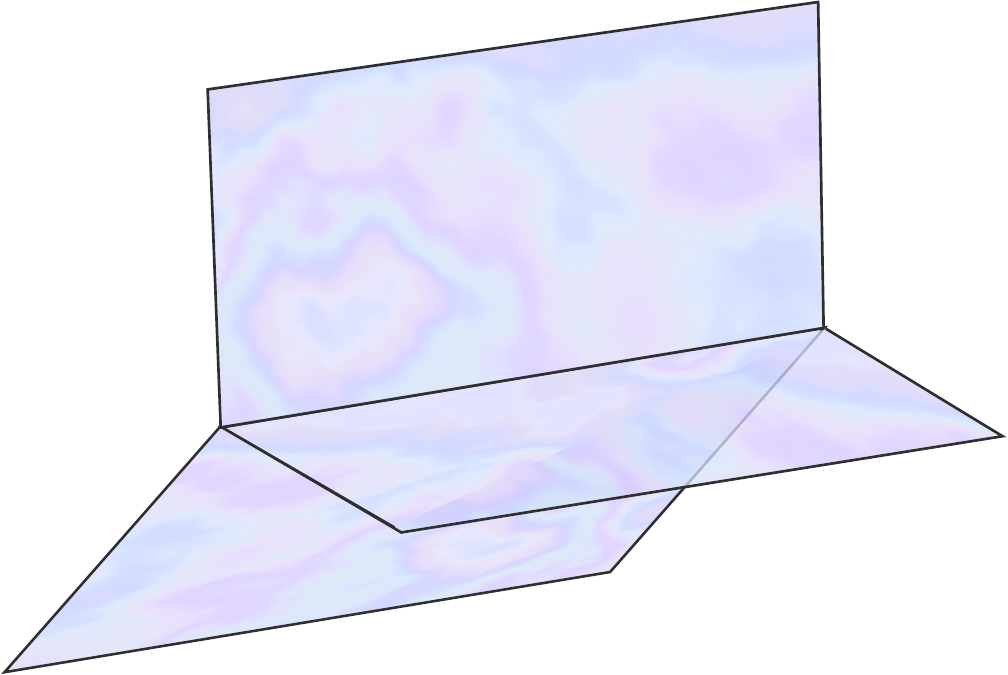}
	\caption{\emph{Three surfaces meeting at a triple junction singularity}}
	\label{fig:TripleJunction}
\end{figure}

\begin{figure}
	\centering
	\subfloat[A wireframe as a boundary]{%
		\includegraphics[scale=0.5]{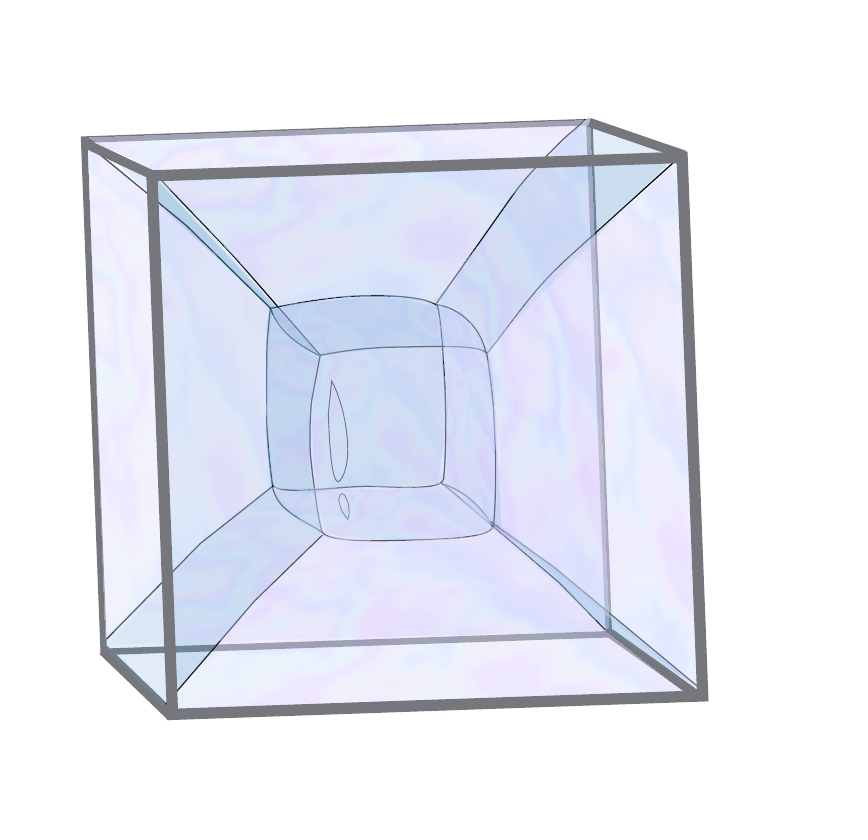}
		\label{fig:cube}}
	\quad
	~
	\subfloat[A non-closed boundary]{%
		\includegraphics[scale=0.15]{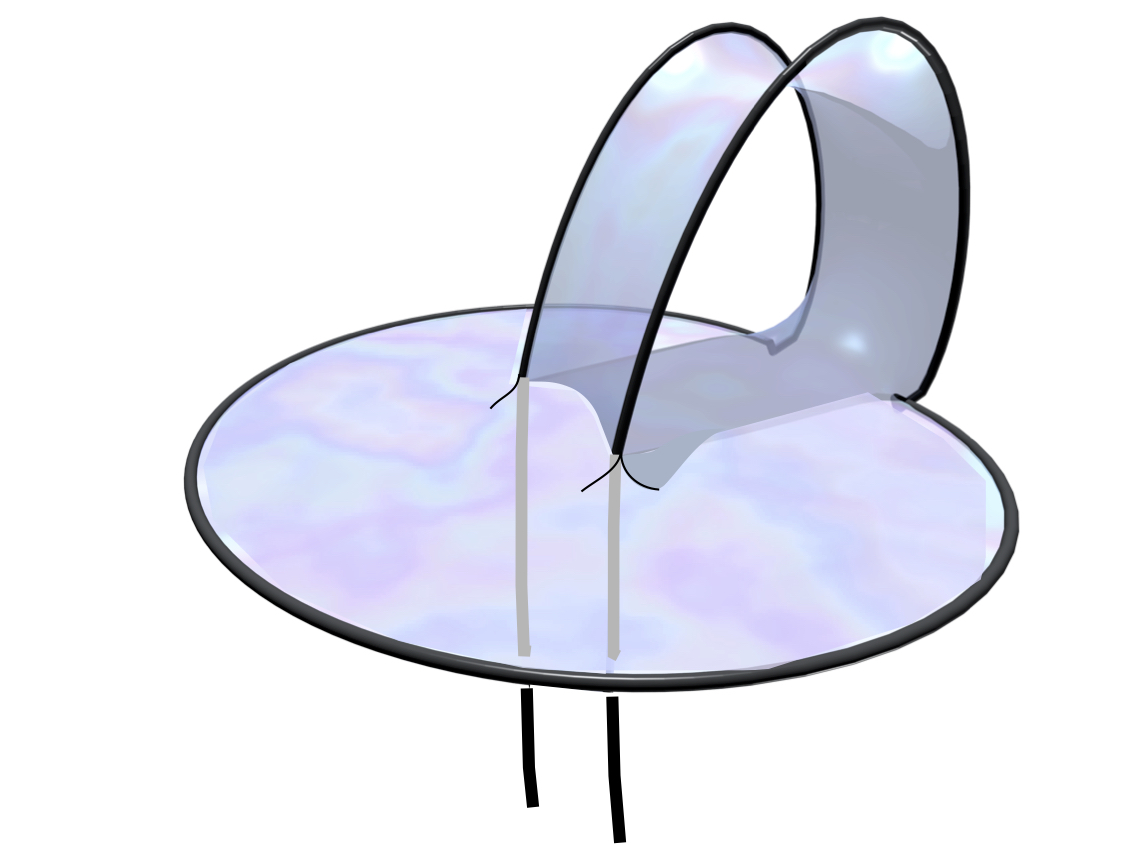}
		\label{fig:Almgren}}
	\caption{\emph{Non-standard boundaries}}
	\label{fig:nonstandardboundaries}
\end{figure}

\begin{figure}
	\centering
	\subfloat[A soap film which retracts onto its boundary]{%
		\includegraphics[scale=0.08]{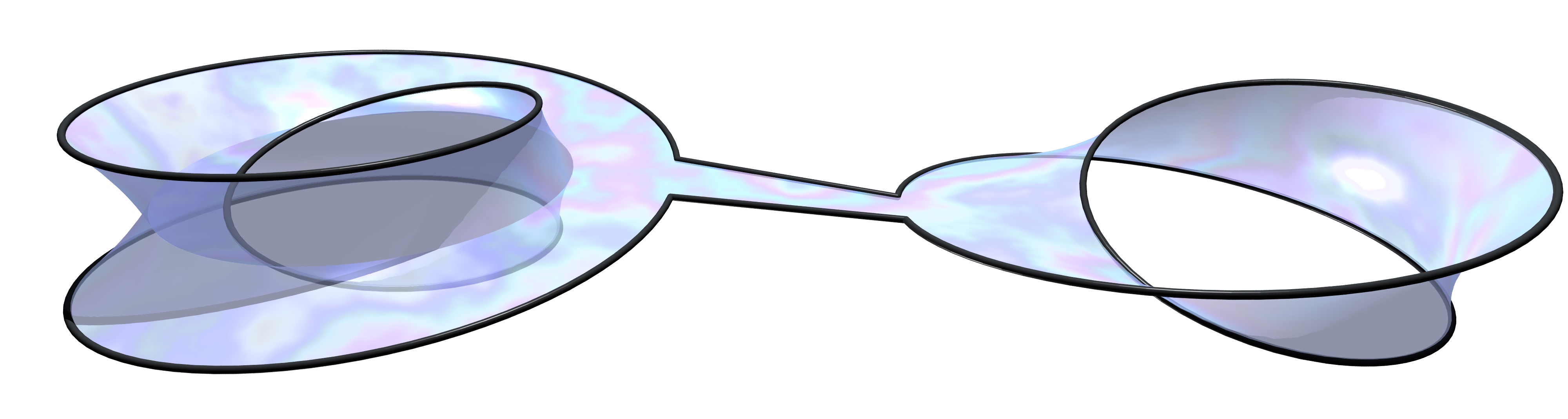}
		 \label{fig:AdamsSurface}}
		~
		
	\subfloat[The boundary wire]{%
		\includegraphics[scale=.08]{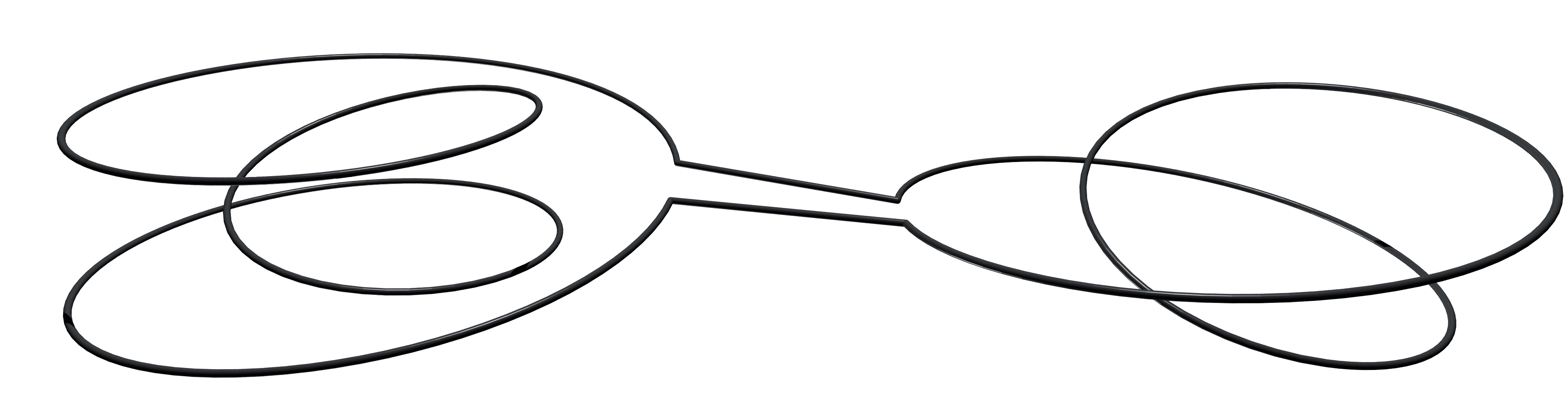}
		\label{fig:exp}}
	 \centering
	\caption{\emph{The Adams surface (A) retracts onto its boundary (B). The left portion of (A) is a triple M\"{o}bius band (a ``Y'' cross an interval, glued along the ends with a \( 1/3 \) twist, see also Figure \ref{fig:TripleMobiusStrip}) and the right portion is a classical M\"{o}bius band. They are joined by a bridge, so that the boundary is a single Jordan curve.}}
	\label{fig:Adams}
\end{figure}

Plateau's problem was first posed by Lagrange, who in 1760 derived the minimal surface equation and asked if one could find a surface of minimal area with a prescribed boundary. The problem was later named after Joseph Plateau \cite{plateauoriginal} who undertook a physical study of soap films and characterized their properties, most notably their singularities. It was actually Lebesgue who coined the term ``Plateau's problem,'' the crux of which was to describe these soap films in mathematical terms. These objects do not behave like classical surfaces such as embedded, or even immersed manifolds with branch points. Three sheets can come together along a line and form what is known as a triple junctions (Figure \ref{fig:TripleJunction}.) Soap films can span wires that are not cycles (Figure \ref{fig:nonstandardboundaries}.) Some films are local area minimizers, yet can retract onto their boundaries (Figure \ref{fig:Adams}.) 

	The long journey towards a full understanding these phenomena began with a much simpler question about the existence of a function with prescribed values on the boundary of a domain \( \O\subset \R^2 \), such that the graph of the function on the interior of \( \O \) is a minimal surface. This problem was studied by Weierstass and Riemann and evolved into the classical theory of minimal surfaces. The next step up in generality came with the study of surfaces defined as images of disks. Jesse Douglas won the first Fields Medal for his solution, which proved the existence of a minimally immersed disk in \( \R^3 \) with a prescribed contour boundary. Many others continued on with striking results of existence and regularity along the way, and the Douglas-Plateau problem for surfaces with higher (non-infinite) genus in arbitrary dimension and codimension was finally solved by Jost in 1985.

\begin{figure}
	\centering
	\subfloat[A minimal surface with infinite genus]{%
		\includegraphics[scale=0.06]{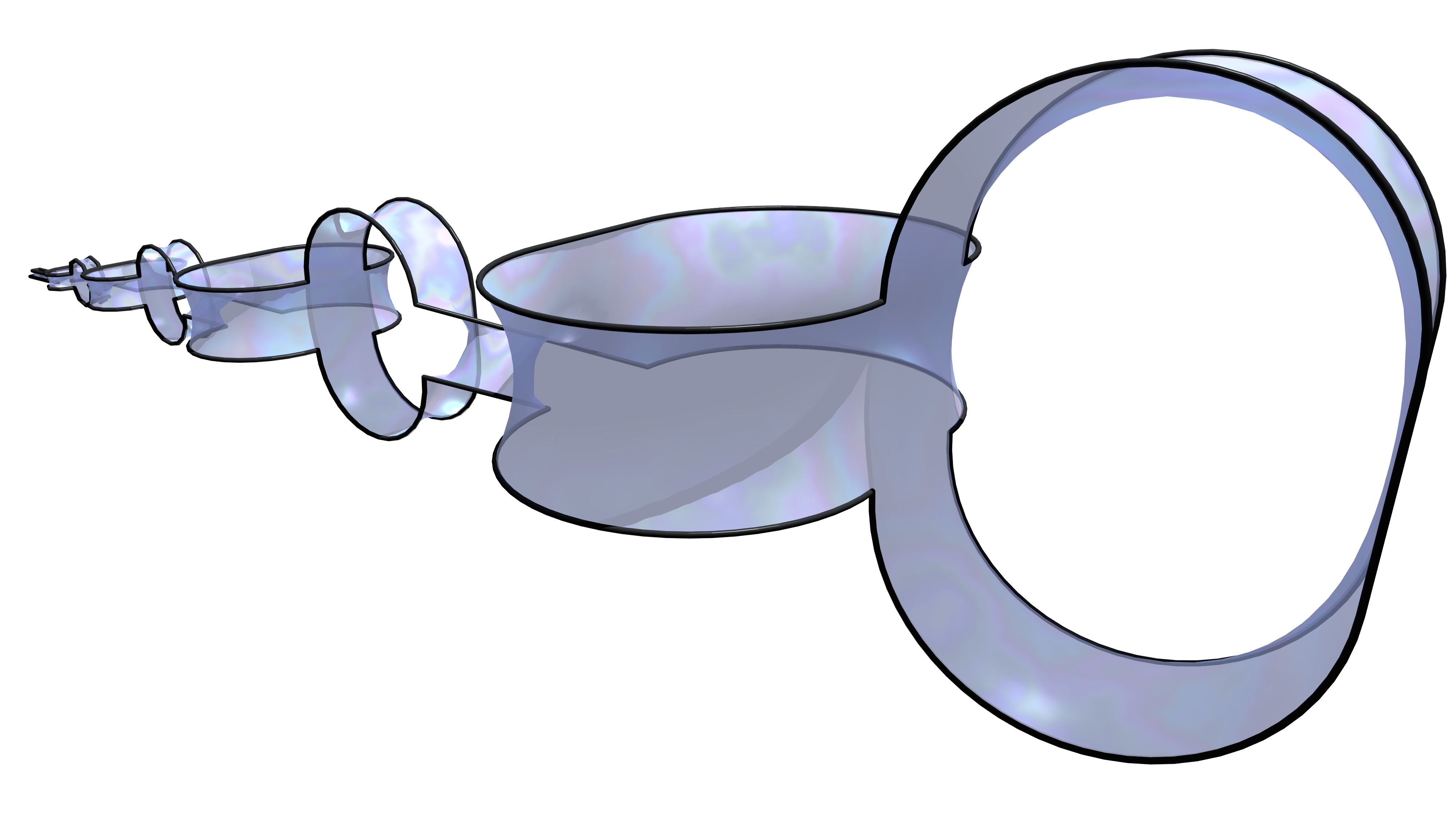}
		\centering
		\label{fig:FlemingSurface}}
		
	\subfloat[The boundary wire can be made smooth except at a single point]{%
		\includegraphics[scale=0.06]{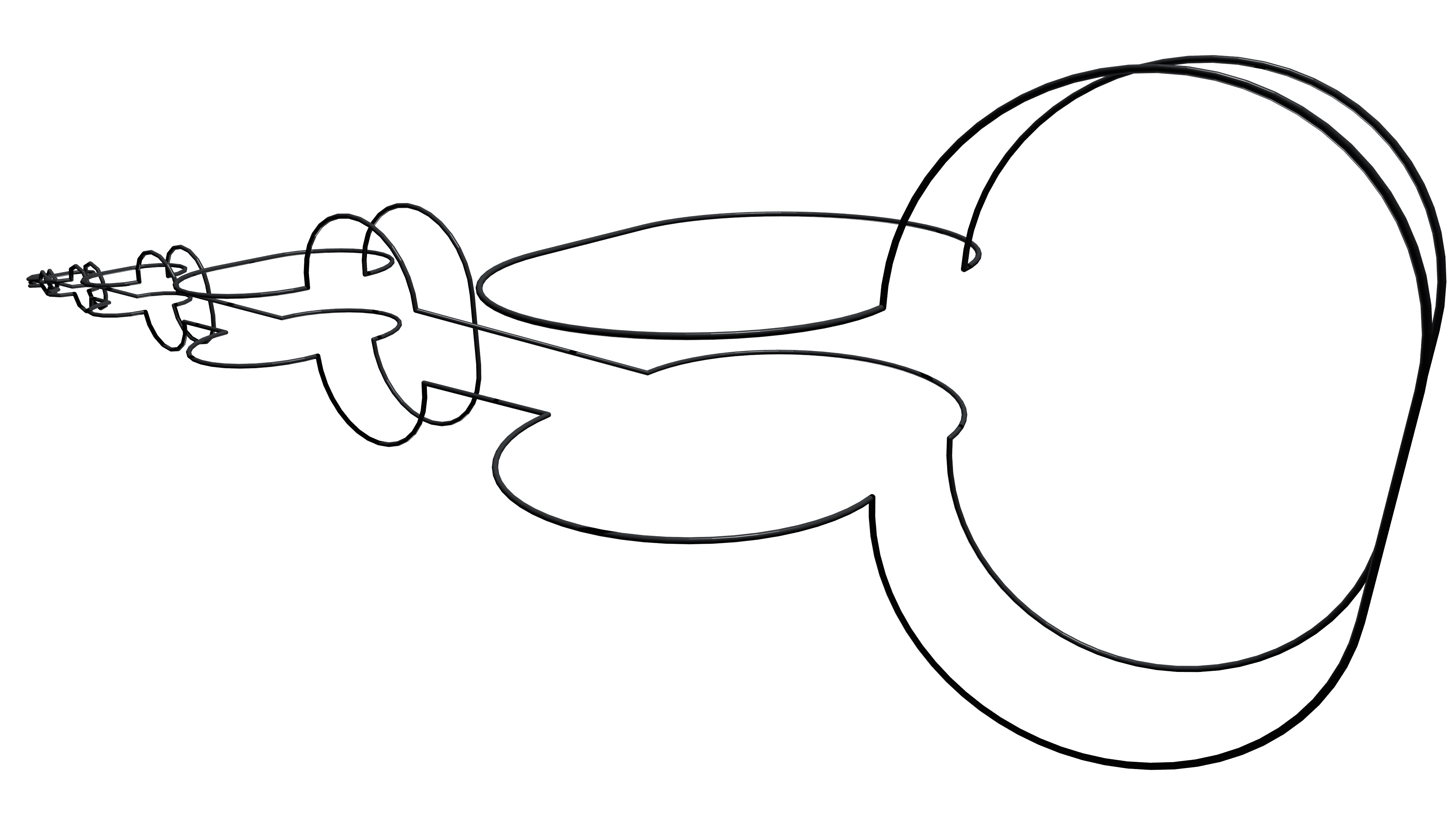}
		\centering
		 \label{fig:FlemingWire}}

	\caption{\emph{Fleming's example of a minimal surface with infinite genus}}
	\label{fig:Fleming}
\end{figure}

However, Fleming demonstrated the existence of a contour boundary which bounds a minimal surface of infinite genus (Figure \ref{fig:Fleming}.) In 1960, Federer and Fleming introduced objects known as integral currents which could model these somewhat pathological surfaces. Their novel approach won them the Steele prize and helped launch the modern study of geometric measure theory. They proved the existence of an integral current with a given boundary which minimizes mass, a quantity which can be thought of as area weighted by an integer multiplicity. It later became known that in low enough dimension, their minimizing current corresponded to an embedded minimal submanifold.

At the same time, Reifenberg thought of a completely different approach. Building on results of Besicovitch and aided by Adams, he defined what it meant for a surface to span a bounding set \( A \) using \v{C}ech homology. Using his theory, he proved the existence of a surface with minimal area amongst those surfaces which can be written as a nested union of manifolds whose boundaries converge to the contour. These surfaces include non-orientable surfaces, as well as the example of Fleming, amongst others. His work is considered to be a masterpiece and deeply influenced several mathematicians, including Morrey, Almgren, Fomenko, as well as ourselves.

Almgren proposed three approaches to Plateau's problem. The first used varifolds \cite{almgrenmimeo, varifolds}, which Young \cite{youngI, youngII} had discovered but called ``generalized surfaces.'' Integral varifolds have a compactness theorem \cite{allard} which can be used to prove the existence of a stationary varifold with smallest area. Integral varifolds model just about any imaginable minimal surface, including the example of Adams (Figure \ref{fig:Adams}.) However, Almgren did not prove that this smallest stationary varifold was the smallest among a class of surfaces which also included non-stationary varifolds.

Almgren's second approach \cite{almgrenannals} was an attempt to generalize Reifenberg's results to elliptic integrands. To read his paper requires  expertise in methods of geometric measure theory, varifolds, integral currents, and flat chains, for it blends them all. It has some gaps, one of which seems serious (see \cite{elliptic}.) 

Almgren's third and final approach \cite{almgren} was to define a new class of surfaces which later became known as quasiminimal sets, which, roughly speaking, have a controllable increase of area under small deformations. Although he was unable to prove an existence theorem of an area minimizer in this category,  his regularity results are of major importance. 

    The authors have recently announced the first solution to the full elliptic Plateau problem  \cite{elliptic}.  Our proof of existence of minimizers builds upon classical measure theory, and techniques of Reifenberg \cite{reifenberg}, Federer and Fleming \cite{federerfleming}.  Our proof of regularity relies upon Almgren \cite{almgren}, although Reifenberg \cite{reifenbergepi, reifenberganalytic} is closely related.  Spanning sets can be defined using homology, cohomology or homotopy. An axiomatic approach without requiring any definition of a spanning set is also provided.  Our results carry over to ambient spaces of Lipschitz neighborhood retracts, including manifolds with boundaries and manifolds with singularities. 

\begin{figure}
	\centering
	\includegraphics[width=.4\textwidth]{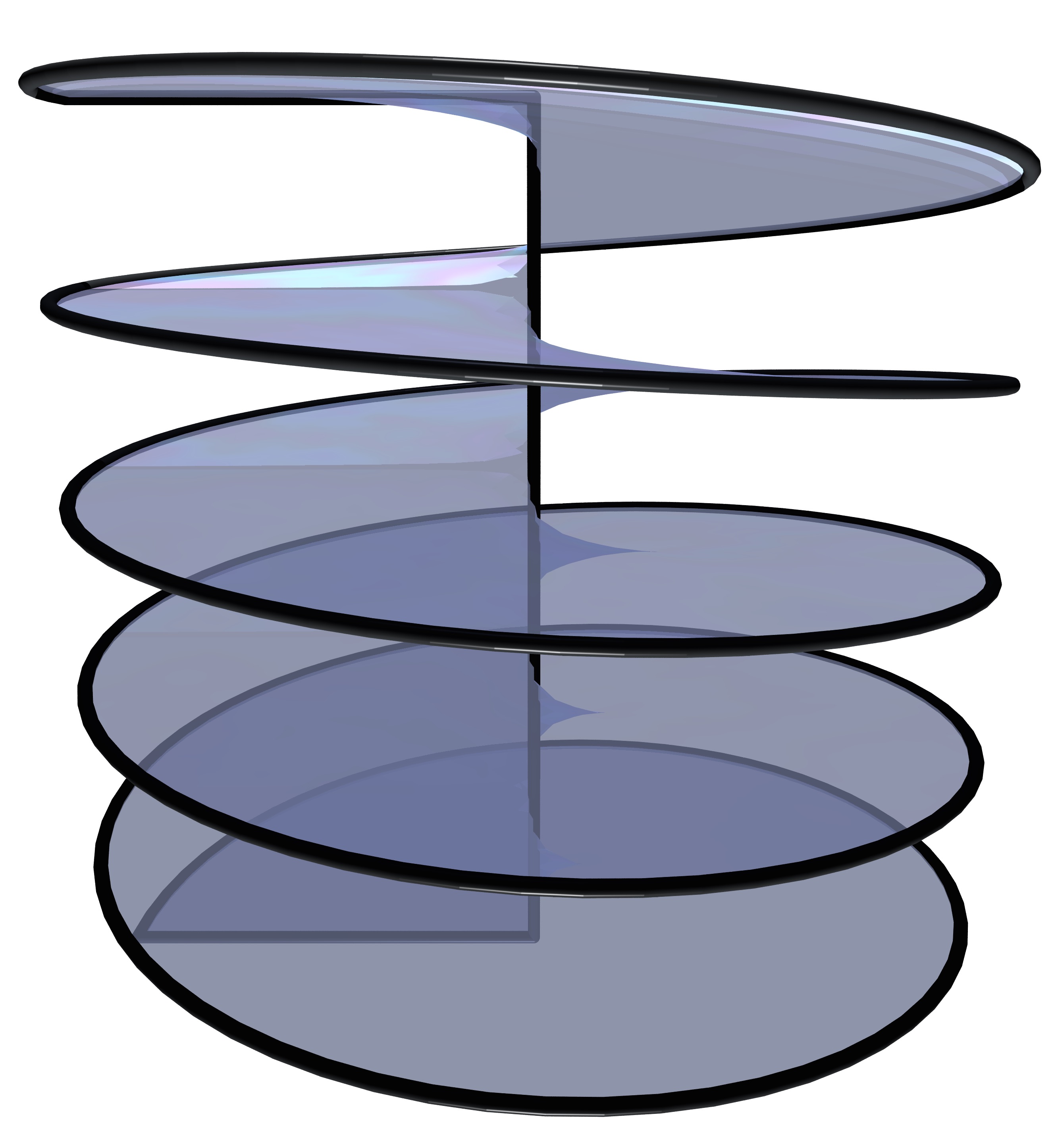}
	\caption{\emph{A piece of the helicoid, which extends to infinity in all directions}}
	\label{fig:Helicoid}
\end{figure}

In the last few years there has been a flurry of other activity in Plateau's problem and in related fields. Papers have recently appeared on sliding boundaries, where the soap film's interface within a larger boundary is permitted to move freely, on flexible boundaries, where the boundary itself is permitted to move subject to forces created by the spanning soap film, and on axiomatic theory, ellipticity, spanning conditions, and classical minimal surface theory. Indeed, very recently, a beautiful and central theorem in classical minimal surfaces was proved by Meeks and Rosenberg, namely that every simply connected, properly embedded minimal surface in \( \R^3 \) must be either a plane or a helicoid (Figure \ref{fig:Helicoid}.)

There are still many major questions left unanswered in Plateau's problem, indeed the list seems to be growing, not shrinking. We have enumerated several of our favorites at the end of this paper, some of which are newly posed.

\subsubsection*{Disclaimer}
\begin{figure}
			\centering
			\subfloat[Two disks]{%
			\includegraphics[width=.4\textwidth]{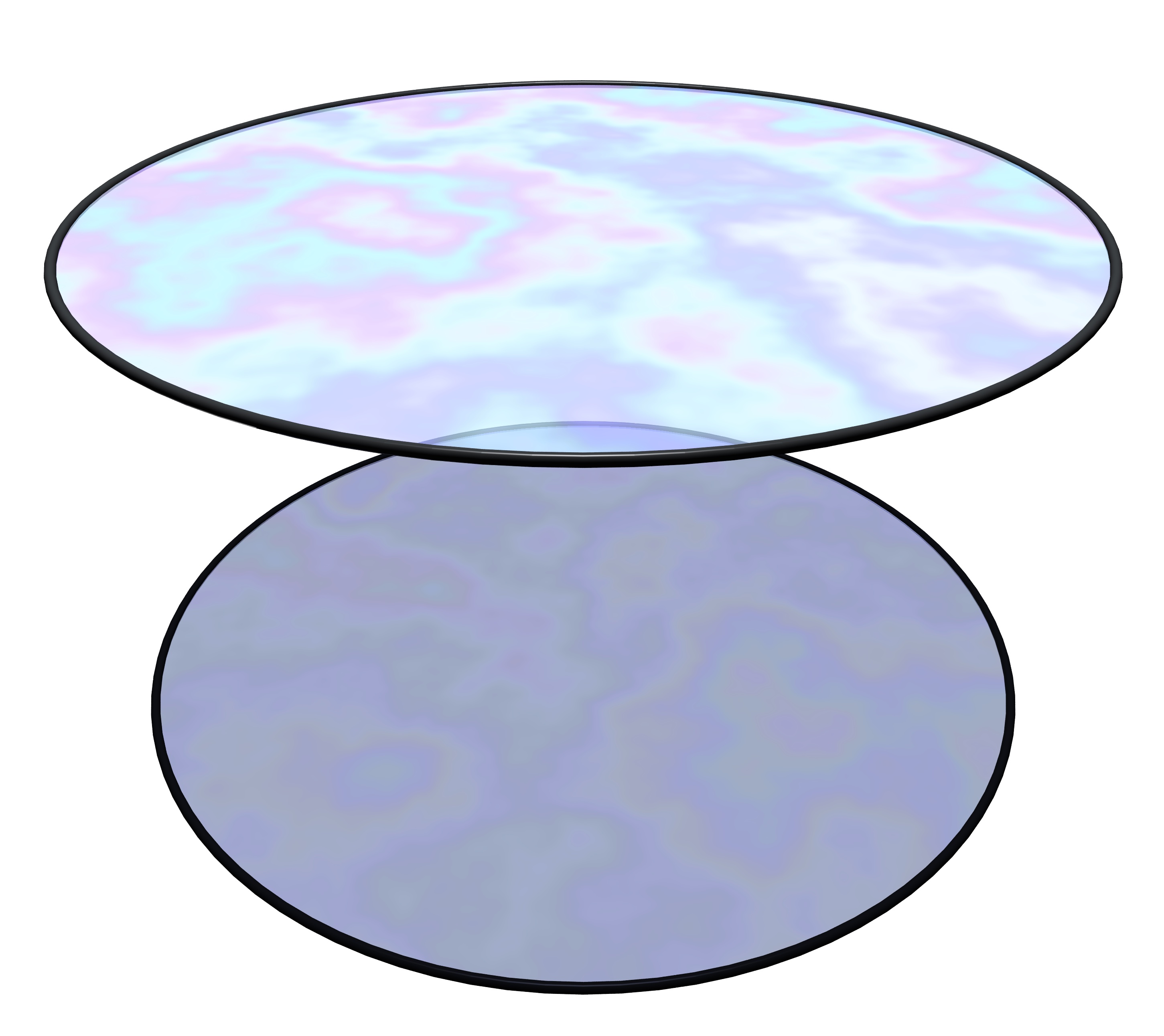}
			\label{fig:TwoDisks}}
			\quad
			~
			\subfloat[Catenoid]{%
			\includegraphics[width=.42\textwidth]{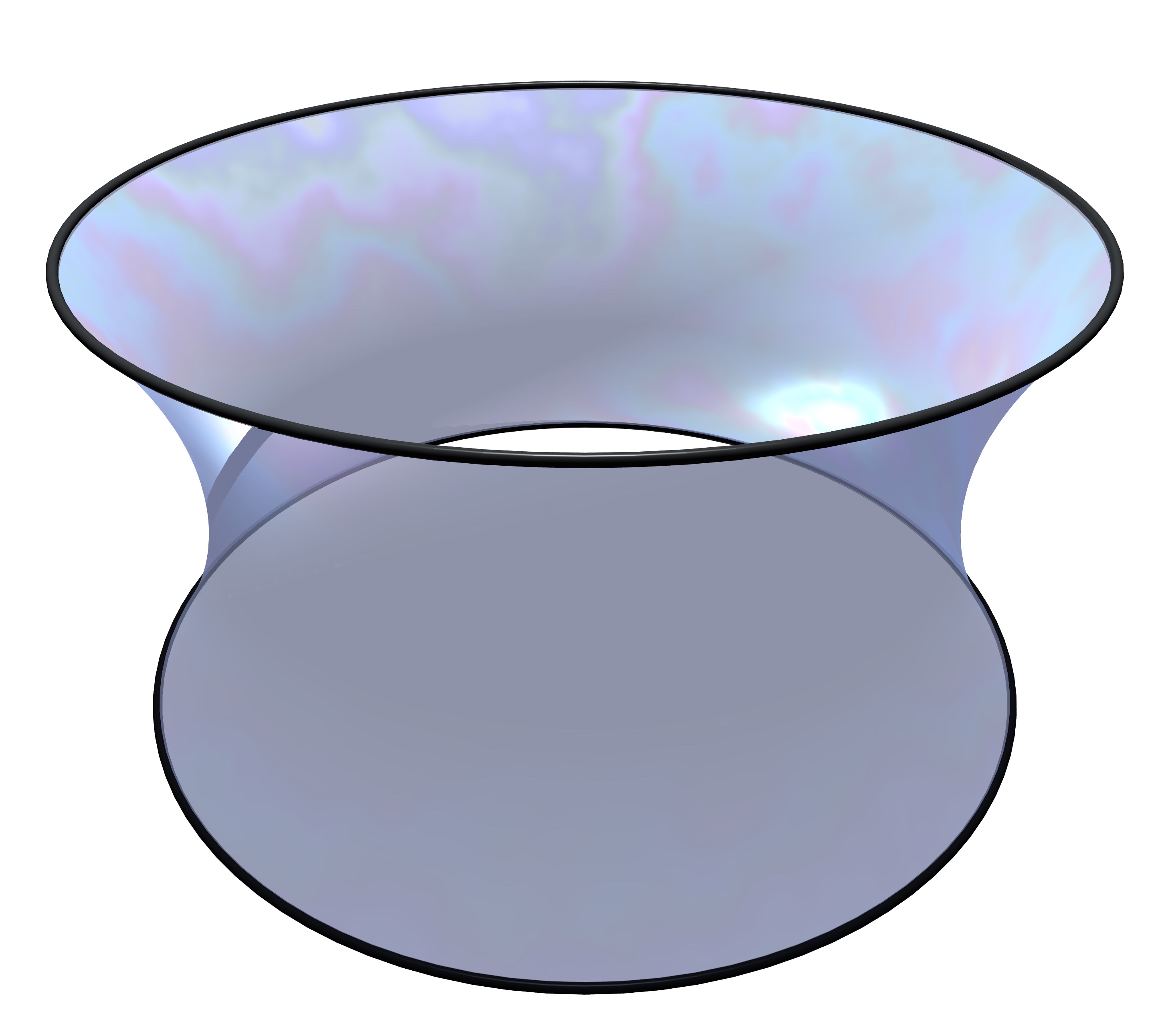}
			\label{fig:Catenoid}}
			\centering
			\caption{\emph{Minimal surfaces with the same boundary}}
		\end{figure}

Before diving in, the reader should be aware that the following exposition is far from complete and will often be imprecise. We hope that it will be useful to give non-specialists some idea of the history of Plateau's problem, a few lines of current theoretical development, and some open problems, both enduring and emerging. We have endeavored to give a broad overview of many different aspects of the problem, and in doing so, have, by necessity, left out many important contributions by numerous mathematicians. If we have neglected to mention your favorite result, it was not due to malice, but rather due to the constraints of writing this article. If we have incorrectly stated your favorite result, know that we are experts in only a small portion of the Plateau problem, and would welcome any corrections you might provide. Lastly, the authors would like to thank Emanuele Paolini\footnote{Figures \ref{fig:Almgren}, \ref{fig:Adams}-\ref{fig:UnorientedProblem}, \ref{fig:TripleMobius}}, Ken Brakke\footnote{Figure \ref{fig:BorromeanTests}}, and Claire-Audrey Bayan\footnote{Figure \ref{fig:TripleJunctionResoltuion}} for the use of their soap-film figures.

\section{Classical Minimal Surfaces}
	\label{sec:minimal_surfaces}	
	\subsection{The Minimal Surface Equation}
		\label{sub:the_early_years}
		Let \( \Sigma \) be a surface in \( \R^3 \). We say that \( \Sigma \) is a \emph{\textbf{minimal surface}} if every point in \( \Sigma \) has an \( \epsilon \)-neighborhood \( U \) which has least area among all surfaces \( S\subset \R^3 \) with boundary \( \p U \). This condition is equivalent to \( \Sigma \) having vanishing mean curvature, and to the condition that \( \frac{\p}{\p t} \textrm{Area}(\Sigma_t) \lfloor_{t=0}=0 \) for all compactly supported variations \( \Sigma_t \) of \( \Sigma \). These are the two (and higher) dimensional analogs of geodesics, but one must be careful in this comparison: even in \( \R^n \), a minimal surface with a given boundary may not have smallest area amongst all surfaces with that boundary. For example, consider two horizontal disks in \( \R^3 \) separated vertically by a small amount. Their union forms a minimal surface, but the cylinder has smaller area, and the catenoid, smaller still.  
				
		A special kind of minimal surface \( \Sigma \) is one which occurs as a graph:

		Suppose \( \O \) is a bounded open set in \( \R^2 \) with locally Lipschitz boundary \( \p \O \), and suppose \( g:\p \O\to \R \) is continuous. Let \( \mathit{F} \) denote the set of continuous extensions of \( g \) to \( \O \) which are continuously differentiable in \( \O \) and whose derivative is integrable. Let \( A(f) \) denote the surface area of the graph of such a function \( f \) and suppose \( f\in \mathit{F} \) solves Plateau's problem for this setup. That is, \( f \) satisfies
		\begin{equation}
			\label{eq:star}
			A(f)\leq A(h) \quad\quad \textrm{for all } h\in \mathit{F}.
		\end{equation}
		Then,
		\begin{equation}
			\label{eq:mse}
			\mathit{div} \left(\frac{\nabla(f)}{\left(1 + |\nabla f|^2 \right)^{1/2}}\right) = 0.
		\end{equation}
		
		This differential equation is called the \emph{\textbf{minimal surface equation}}. It is the Euler-Lagrange equation for the area functional.

		\begin{thm}
			\label{thm:dir}
			A function \( f \in \mathit{F} \) satisfies \eqref{eq:star} if and only if \( f \) satisfies the minimal surface equation \eqref{eq:mse}.
		\end{thm}
		
		Such a function \( f \) is unique, and is in fact analytic on \( \O \) \cite{hopf}.

		Suppose now the domain \( \O \) is a disk \( D \) of radius \( r \). This produces a unique minimal surface for each continuous function \( g \) on \( S^1 \), and thus there are uncountably many minimal surfaces that are graphs over disks of radius \( r < \i \). However, if \( r = \i \), the situation simplifies dramatically:
		\begin{thm}
			\label{thm:Bernstein}(Bernstein \cite{bernstein})
			Any solution of the minimal surface equation \eqref{eq:mse} which is defined on all \( \R^2 \) must be linear. 
		\end{thm}
		Bernstein conjectured that this was also the case in higher dimension. Indeed this is the case up to dimension seven:
		\begin{thm}
			\label{thm:adg}(de Giorgi \cite{degiorgibernstein}, Almgren \cite{almgreninterior})
			The Bernstein conjecture holds for minimal graphs \( \G = \{(x,f(x)) \in \R^{n+1} :x \in \R^n\} \) for \( n \leq 4 \).
		\end{thm}
		\begin{thm}
			\label{thm:simons}(Simons \cite{simonsminimal})
			The Bernstein Conjecture holds for minimal graphs when \( n \leq 7 \).
		\end{thm}
		
	\subsection{Recent Developments}
		\label{sub:recent_developments}
		In the classical theory, non-compact minimal surfaces might not be graphs and might not have a boundary. An example in \( \R^3 \) is the helicoid (Figure \eqref{fig:Helicoid}). A long outstanding question posed by Osserman was the following generalization of Bernstein's Conjecture: The plane and the helicoid are the only properly embedded, simply-connected, minimal surfaces in \( \R^3 \).

		Osserman's conjecture has recently been solved after many years of effort by Meeks and Rosenberg \cite{meeksrosenberg} who built on the work of Colding and Minicozzi \cite{coldingminicozzi}, \cite{coldingminicozziiv} as well as a number of other mathematicians (see \cite{meeksperez} for a detailed discussion and a more complete list of citations.)

		Combining this with work of Collin \cite{collin}, L\'opez and Ros \cite{lopezros}, Meeks, P\'erez and Ros \cite{meeksperezros} give the following classification theorem:
		\begin{thm}
			\label{thm:catenoid}
			Up to scaling and rigid motion, any connected, properly embedded, minimal surface in \( \R^3 \) is a plane, a helicoid, a catenoid or one of the Riemann minimal examples. In particular, for every such surface there exists a foliation of \( \R^3 \) by parallel planes, each of which intersects the surface transversely in a connected curve which is a circle or a line.
		\end{thm}

		This beautiful theorem is only just the beginning of what looks to be a new era in classical minimal surface theory. A few particularly intriguing open problems are listed in  \S\ref{sec:openproblems} and can be found as part of a larger list in \cite{meeksperez}.
		
\section{The Douglas-Plateau Problem: Immersions of Disks and Surfaces of Higher Genus}
	\label{sec:douglas}

	The Plateau problem, in its original formulation, was to find a minimal immersed disk whose boundary was a given Jordan curve in \( \R^n \). For any immersed disk, coordinates on the disk \( D \) can be chosen so that the immersion is conformal, in which case minimality is equivalent to the immersion being harmonic.

	Independently, Jesse Douglas \cite{douglas} and Tibor Rad\'o \cite{rado} proved the existence of such a surface, with (possibly) isolated singularities:
	
	\begin{thm}
		\label{thm:douglas}
		If \( C \) is a Jordan curve in \( \R^n \), there exists a continuous map \( \iota: D\to \R^n \) which is conformal and harmonic away from a set of isolated singularities (i.e. branch points,) such that \( \iota\lfloor_{\p D} \) parameterizes \( C \).
	\end{thm}

	The solutions produced by Douglas and Rad\'o also had minimal area in the class of branched immersions. However, Douglas could prove slightly more than Rad\'o, principally his theorem allowed for certain pathological boundaries \( C \) which could only be spanned by disks of infinite area\footnote{Douglas was well known for his displeasure at having to share credit with Rad\'o for Theorem \ref{thm:douglas}. When teaching subsequent geometry courses, he eschewed those books which, in covering the theorem, contained the attribution ``Douglas-Rad\'o.'' Unfortunately, as the years went on, he was forced to use increasingly antiquated texts, since virtually no book published after the 1930s failed to give this (correct) attribution (See \cite{rassias}. The second author has also heard a similar story from Martin Bendersky, who was a student in one of Douglas's courses at City College.)}. In addition, Douglas's methods signaled a significant and promising departure from the classical techniques, and for these reasons he was awarded the first Fields medal in 1936.

	Osserman proved \cite{osserman} that if \( n=3 \), then branch points did not exist in such minimal disks. Thus,
	
	\begin{thm}
		If \( C \) is a smooth Jordan curve in \( \R^3 \), there exists a (conformal, harmonic) immersion \( \iota: D\to \R^3 \) such that \( \iota\lfloor_{\p D} \) parameterizes \( C \), whose area is minimal among all immersed disks whose boundaries parameterize \( C \).
	\end{thm}
	
	\begin{figure}
		\centering
		\subfloat[An immersed disk spanning the boundary of the M\"{o}bius band]{%
		\includegraphics[width=.47\textwidth]{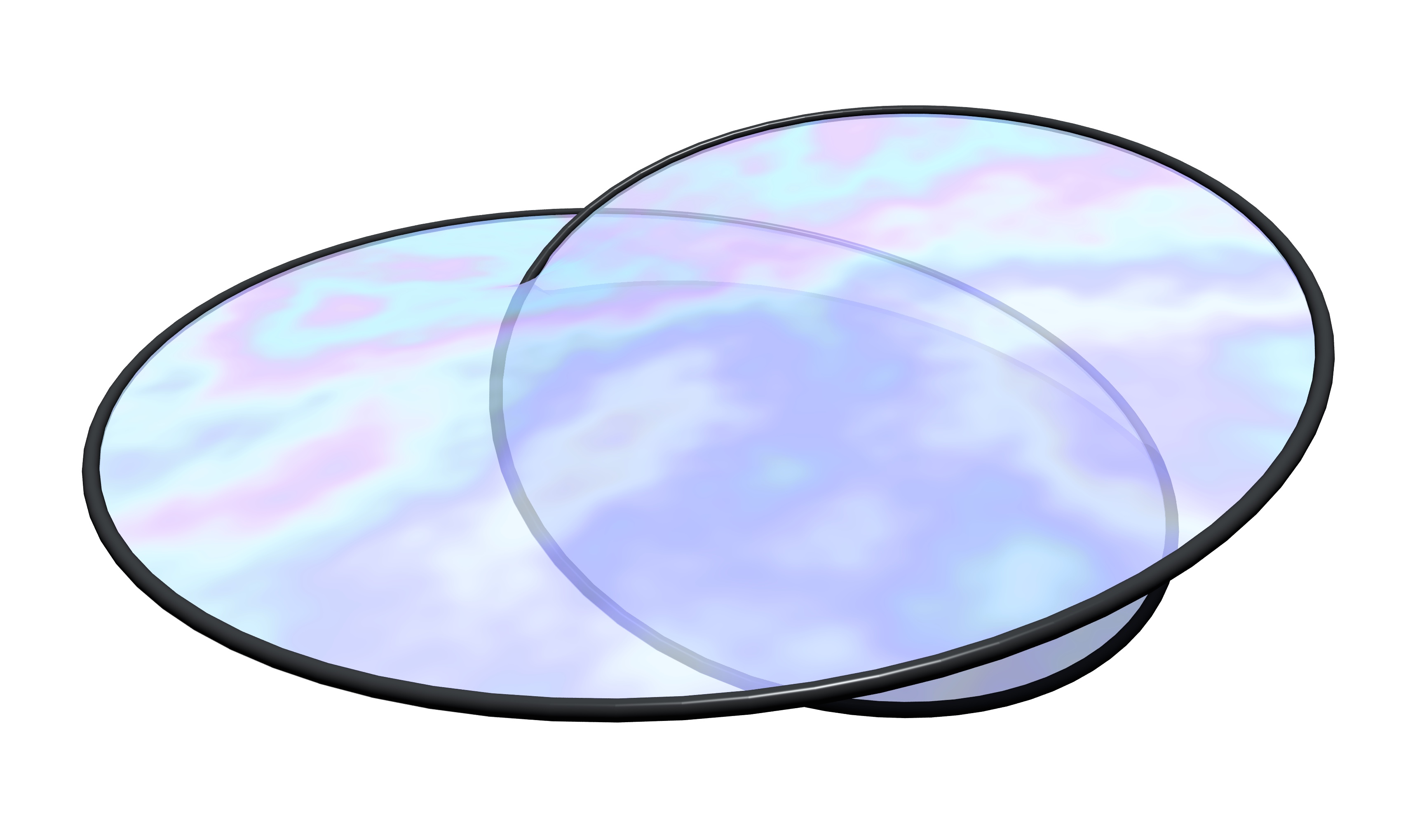}
		\centering
		\label{fig:OrientedMobiusSpanner}}
		\quad
		~
		\subfloat[A M\"{o}bius band]{%
		\includegraphics[width=.47\textwidth]{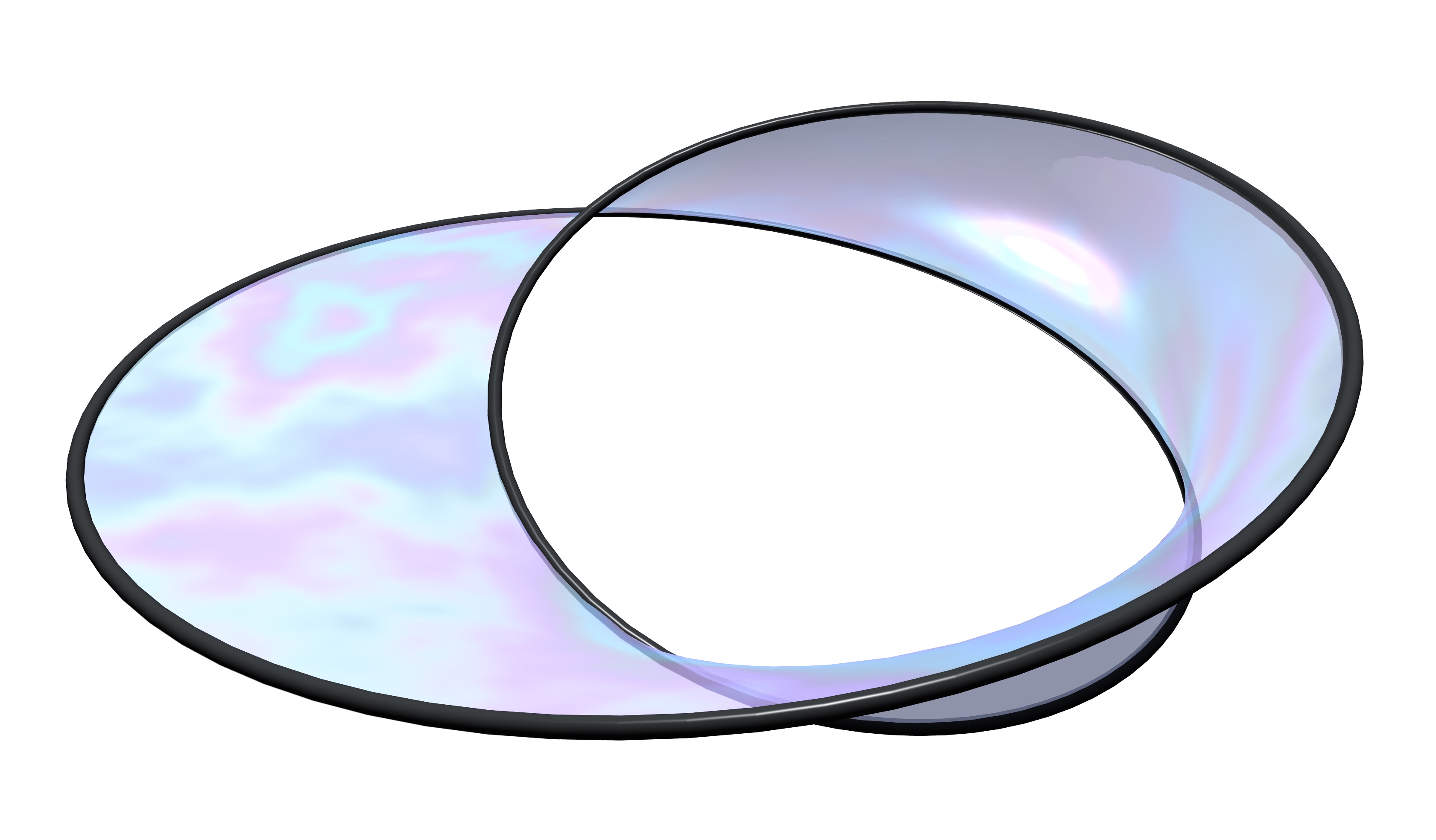}
		\centering
		\label{fig:MobiusStrip}}
		\caption{\emph{A non-orientable surface with smaller area than any immersed disk with the same boundary}}
		\label{fig:UnorientedProblem}
	\end{figure}
	
 	However, there still could be immersed surfaces with boundary \( C \) whose area is strictly less than those of the solutions produced above. Consider a thin M\"{o}bius band. Its area is less than that of any immersed disk spanning the boundary curve (Figure \ref{fig:UnorientedProblem}.) Higher genus surfaces could also have less area. Douglas had made attempts to generalize his techniques to account for possibly non-orientable surfaces and those of higher genus, but the consensus seems to be that his arguments were incomplete \cite{gray}. It took until the 1980s for a complete solution to the higher genus problem to appear in a paper by Jost \cite{jost} (see \cite{bernatzki} for the non-orientable case,) who built upon ideas of Schoen-Yau \cite{yau} and Sacks-Uhlenbeck \cite{uhlenbeck}. Tomi-Tromba \cite{tromba} soon after offered a different solution to the higher genus problem based on a development of Teichm\"{u}ller theory from the viewpoint of differential geometry. 

	\begin{figure}
	  \centering
	    \includegraphics[width=.85\textwidth]{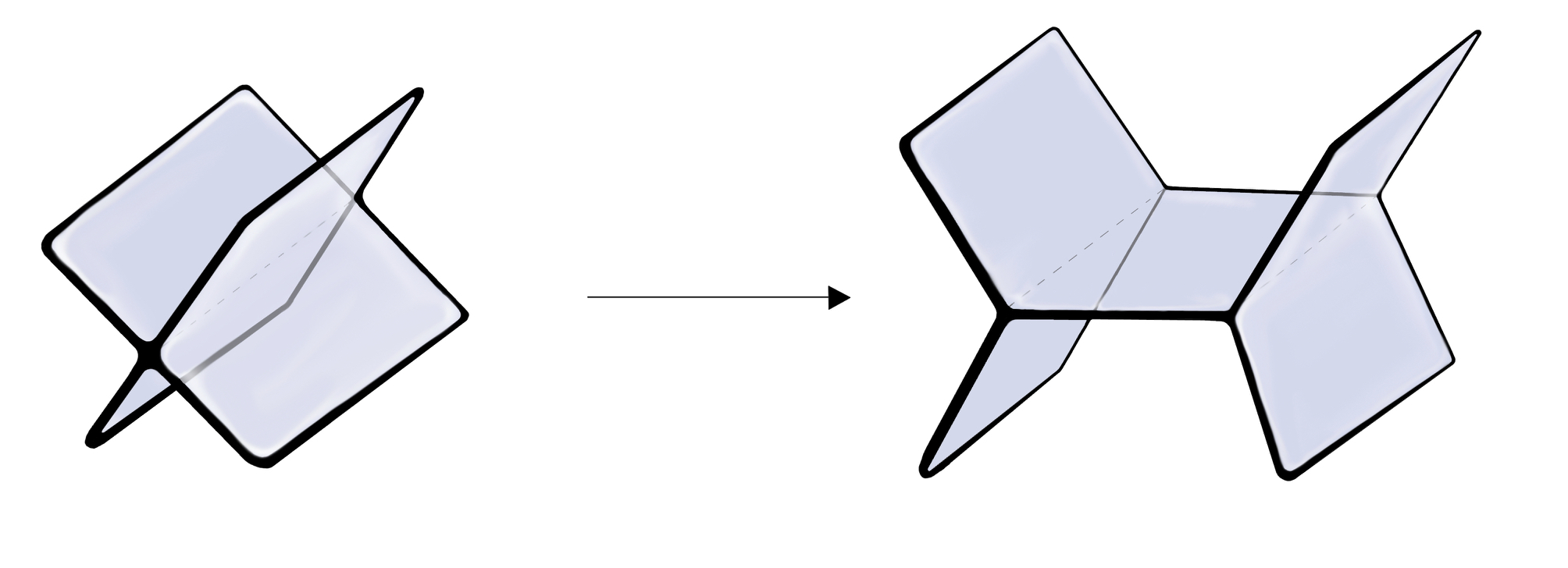}
	\caption{\emph{The transverse intersection of two sheets can be replaced with a pair of triple junctions with smaller area. These singularities do not, however, show up in the mass minimization problem (see \S\ref{sub:massminimization},) as the horizontal portion of the right hand figure would require higher mass to cancel out any contribution of the triple junction to the boundary.}}
	\label{fig:TripleJunctionResoltuion}
	\end{figure}
	
	Although not a mathematical shortcoming, self-intersections in which two sheets intersect transversally can easily show up in Douglas's and others' immersed solutions. Such solutions are not physically realistic as soap-films, since transverse intersections resolve into pairs of triple junctions (Figure \ref{fig:TripleJunctionResoltuion}) with smaller total area when one is allowed to consider surfaces more general than immersed manifolds. These generalized surfaces with triple junctions became important in the 1960s in the work of Reifenberg \cite{reifenberg}, Almgren \cite{almgrenannals} and Taylor \cite{taylor}, who studied soap-film regularity and classified the singularities for these general size-minimizing surfaces.

\section{Orientable Generalized Surfaces}
	Fleming's example (Figure \ref{fig:Fleming}) shows that one must consider all topological types to find a true minimizer for an orientable version of Plateau's problem. The curve in the figure is an unknotted simple closed curve and is smoothly embedded except at one point. It bounds the shaded surface which is orientable and clearly area minimizing. This shows that one should not insist that competitors have finite topological type when looking for absolute area minimizers. The proof of the existence of an orientable surface minimizing the area among all possible surfaces without restriction on their topological type requires other techniques, namely the results of Reifenberg  \cite{reifenberg} and the integral currents of Federer and Fleming \cite{federerfleming}.
			
	\subsection{Hausdorff Measure}

		Readers will recall that for any \( E \subset \R^n \), and any non-negative real number \( k \), the \( k \)-dimensional Hausdorff measure of \( E \) is \[ \H^k(E) = \lim_{\d \to 0} \H_\d^k(E) \] where \[ \H_\d^k(E) =  \alpha_k \inf \left\{\sum_{i \in I} \textrm{diam}(U_i)^k: E \subset \cup_{i \in I} U_i \right\}, \] \( \alpha_k \) is a normalizing constant, and the infimum is taken over all coverings of \( E \) by a collection \( \{U_i \}_{i\in I} \) of sets with \( \textrm{diam}(U_i) < \d \). The normalizing constant \( \alpha_k \) is chosen so that when \( k \) is an integer, the \( k \)-dimensional Hausdorff measure of the unit cube in \( \R^k \) is one.

		Note that \( \H_\d^k(E) \) is monotone decreasing in \( \d \), so the limit \( \lim_{\d \to 0} \H_\d^k(E) \) exists but may be infinite. \( \H^k \) is a Borel regular outer measure and coincides with Lebesgue measure when \( E \) is a \( k \)-dimensional submanifold.
	
		The \emph{\textbf{Hausdorff dimension}} of \( E \) is the infimum over all \( k\geq 0 \) such that \( \H^k(E)=0 \).
		
	\subsection{Integral Currents and Mass Minimization}
		\label{sub:massminimization}
		A \emph{\textbf{\( k \)-dimensional current}} \( T \) on an \( n \)-dimensional smooth manifold \( M \) is a linear functional on the space of compactly supported smooth \( k \)-forms \( \cal{D}^k(M) \), continuous in the following sense:

		If \( \o_i \) is a sequence of \( k \)-forms supported in a single compact set \( K \) contained in a coordinate neighborhood and \( \p^r \o_i \to 0 \) uniformly for all \( 0\leq |r|<\i \), then \( T(\o_i)\to 0 \). Here \( r=(r_1,\dots,r_n) \) is a \( n \)-tuple of non-negative integers, \( |r|=\sum r_i \), and \( \p^r \) is shorthand for the coordinate-wise differentiation operator \[ \frac{\p^{|r|}}{\p x_1^{r_1}\dots \p x_n^{r_n}}. \]

		It is important to note that the topology on \( \cal{D}^k(M) \) is strictly finer than the subspace topology induced by the inclusion of \( \cal{D}^k(M) \) into the space \( \cal{E}^k(M) \) of \( C^\i \) differential \( k \)-forms on \( M \), in which a sequence of forms \( \eta_i \) converges to zero whenever \( \p^r \eta_i\lfloor_K\to 0 \) uniformly for all \( 0\leq r<i \) and all compact sets \( K \) contained in a coordinate neighborhood. The difference is subtle, but has extremely important consequences (the full description of the space \( \cal{D}^0 \), its topology, and continuous dual was the major component of L. Schwartz's Fields medal.) For example, in \( \cal{E}^0(\R) \), any sequence of bump functions \( f_i \) equal to \( 1 \) on the interval \( [-i,i] \) converges to the function \( 1 \), which is no longer compactly supported. Such a sequence is not convergent in \( \cal{D}^0(\R) \). As a matter of fact, \( \cal{D}^k(\R^n) \) is complete, so it is not even Cauchy.
		
		\subsubsection{Examples}
			\begin{itemize}
				\item If \( S\subset M \) is an oriented \( k \)-dimensional submanifold, then a \( k \)-current \( [[S]] \) is defined, setting \( [[S]](\o)\equiv \int_S \o \).
				\item If \( M \) is equipped with a volume form \( dV \), then a \( k \)-vector field \( X \) on \( M \) defines a current \( [[X]] \), whereby \( [[X]](\o)\equiv \int_M \o(X) dV \).
				\item A generalized Dirac delta is a current: if \( p\in M \) and \( \alpha\in \L^k(T_p M) \), then a current \( [[(p,\alpha)]] \) is defined, where \( [[(p,\alpha)]](\o)\equiv \o_p(\alpha) \).
				\item An \( (n-k) \)-form \( \eta \) in \( \cal{E}^{n-k}(M) \) defines a \( k \)-current \( [[\eta]] \), where \( [[\eta]](\o)=\int_M \eta\wedge \o \). Such a current is called a \emph{\textbf{smooth current}}. Through convolution, it is possible to construct a smoothing operator which approximates any current by a smooth current (See \S 15 \cite{derham1}.)
			\end{itemize}

		Denote the space of \( k \)-dimensional currents by \( \cal{D}_k(M) \). The operator dual to exterior differentiation on forms, denoted \( \p \), turns \( \cal{D}_\bullet(M) \) into a chain complex. The image and kernel of \( \p \) are closed. When given the opposite grading (i.e. give \( \cal{D}_k(M) \) degree \( n-k \),) the resulting cochain complex \( (\cal{D}_{n-\bullet},\p) \) is quasi-isomorphic, via application of the the aforementioned smoothing operator, to the cochain complex \( (\cal{E}^\bullet(M),d) \). Thus, Poincar\'e duality holds in this setting: the homology of currents in degree \( k \) (which is dual to the compactly supported de Rham cohomology in degree \( k \)) is isomorphic to de Rham cohomology in degree \( n-k \).
		
		Before we can describe a sub-complex of \( (\cal{D}_\bullet,\p) \) which computes the integral homology of \( M \), it will be necessary to define the mass of a current. If \( M \) is equipped with a riemannian metric and \( W\subset M \), define \( \|T\|(W)\equiv \sup \{ T(\o): \supp(\o)\subset W,\,\, \|\o\|_0\leq 1 \} \), where \( \|\o\|_0 \) is the supremum of \( \o_p \alpha \), where \( p\in M \) and \( \alpha \) is a unit simple \( k \)-vectors in \( \L^k T_p M \). The (possibly infinite) \emph{\textbf{mass}} of \( T \), denoted \( \mathbf{M}(T) \), is the quantity \( \|T\|(\R^n) \). If \( \|T\|(W) \) is finite for every \( W\subset\subset M \), we say \( T \) has \emph{\textbf{locally finite mass}}.

		A \( k \)-dimensional current \( T \) is called \emph{\textbf{(integer) rectifiable}} if it has locally finite mass and there exists a sequence \( S_i \) of \( C^1 \) oriented \( k \)-dimensional submanifolds of \( M \), a sequence of pairwise disjoint closed subsets \( K_i\subset S_i \) and a sequence of positive integers \( k_i \) such that \[ T(\o)=\sum_i k_i\int_{K_i}\o \] for all \( \o\in \cal{D}^k(M) \). If \( T \) and \( \p T \) are rectifiable, we say that \( T \) is an \emph{\textbf{integral current}}. One can show, e.g. using sheaf theory, that the homology of the chain complex \( (\cal{I}_\bullet(M),\p) \) of integral currents computes the homology of \( M \) with \( \Z \) coefficients.
	
		One can also show that the mass \( \mathbf{M}(T) \) of an integral current \( T \) is the same as the quantity \( \sum_i k_i \H^k(K_i) \).
		
		Central to the utility of integral currents is the following compactness theorem:

		\begin{thm}(Federer-Fleming)
			If \( \{T_i\} \subset \cal{I}_k(M) \) is a sequence of integral currents such that \[ \sup_i \|T_i\|(W)+\|\p T_i\|(W)<\i \] for all \( W\subset\subset M \), then there exists an integral current \( T \) and a subsequence of \( \{T_i\} \) which converges weakly to \( T \).
		\end{thm}

		Since mass is weakly lower-semicontinuous, Federer-Fleming produced the following corollary:
		\begin{cor}
			If \( T\in \cal{I}_k(M) \), then there exists \( T_0\in \cal{I}_k(M) \) with \( T-T_0=\p R_0 \) for some \( R_0\in \cal{I}_{k+1} \) such that \[ \mathbf{M}(T_0)=\inf_{R\in \cal{I}_{k+1}} \mathbf{M}(T_0+\p R). \]
		\end{cor}

		As a special case, setting \( Q=\p T \):
		\begin{cor}
			If \( Q\in \cal{I}_{k-1}(\R^n) \), there exists \( T_0\in \cal{I}_k(\R^n) \) with \( \p T_0=Q \) such that \[ \mathbf{M}(T_0)=\inf_{T\in \cal{I}_k, \p T=Q} \mathbf{M}(T). \]
		\end{cor}

		Another special case occurs when \( T \) is a cycle:
		\begin{cor}
			Each class in \( H_k((\cal{I}_\bullet(M),\p)) \) contains a representative of least mass.
		\end{cor}

	\subsection{Regularity}
		The \emph{\textbf{support}} \( \supp(T) \) of a current \( T\in \cal{D}_k \) is the complement of the largest open set \( U \) for which \( \supp(\o)\subset U \Rightarrow T(\o)=0 \). We say \( p\in \supp(T)\setminus \supp(\p T) \) is an \emph{\textbf{interior regular point}} if there exists \( \e>0 \), a positive integer \( \kappa \) and an oriented \( k \)-dimensional smooth submanifold \( S \) such that \( T(\o)=\kappa [[S]](\o) \) for all forms \( \o \) supported in the ball of radius \( \e \) about \( p \). The remaining points in \( \supp(T)\setminus \supp(\p T) \) are called \emph{\textbf{interior singular points}}, the set of which will be denoted \( \cal{S}(T) \).

		\begin{thm}[Complete Interior Regularity]
			\label{thm:fle}
			If \( T_0\in \cal{I}_n(\R^{n+1}) \), where \( 2\leq n\leq 6 \), and the mass of \( T_0 \) is minimal among all integral currents with the same boundary, then \( \supp(T_0)\setminus \supp(\p T_0) \) is an embedded minimal hypersurface in \( \R^n\setminus \supp(\p T_0) \), and \( \cal{S}(T_0) \) is empty.
		\end{thm}
		
		In 1962, Fleming proved the result for \( n=2 \), so other than regularity at the boundary which was to take another 17 years \cite{hardtsimon}, this result completed the solution of the oriented Plateau Problem in \( \R^3 \) for surfaces of all topological types.

		Almgren \cite{almgreninterior} extended Fleming's theorem to \( n=3 \), and Simons extended it up to \( n=6 \) in \cite{simonsminimal}.

		Also in \cite{simonsminimal} Simons constructed an example which showed that singularities could in fact occur in dimension \( 7 \) and higher. The ``Simons cone'' is the cone over \( S^3 \times S^3 \subset S^7\subset \R^8 \). He showed it was locally mass minimizing, yet has an isolated interior singularity. Immediately after Simons published his example, Bombieri, de Giorgi, and Giusti \cite{bombieridegiorgigiusti} showed in a marathon three-day session\footnote{This story was recently communicated by Simons to the second author.} that \( S \) is in fact globally mass minimizing. As a corollary, they also showed that, for any \( n \geq 8 \), there exist functions which satisfy the minimal surface equation and are not affine, finally settling the Bernstein problem in all dimensions.

		Not long after, Federer \cite{federerregularity} put a bound on the size of the singular set \( \cal{S}(T) \):

		\begin{thm}
			\label{thm:fed}
			The singular set \( \cal{S}(T_0) \) has Hausdorff dimension at most \( n - 7 \). Singularities are isolated points if \( n = 7 \).
		\end{thm}

		Bombieri, de Giorgi, and Giusti \cite{bombieridegiorgigiusti} showed that this bound is sharp: there exist mass minimizers \( T_0 \) in every dimension \( n\geq 7 \) such that \( \cal{H}^{n-7}(\cal{S}(T_0))>0 \). In the 90's, Leon Simon \cite{simon95} proved this singular set is well-behaved:

		\begin{thm}
			\label{thm:sim95}
			Except for a set of \( \H^{m-7} \)-measure zero, the singularity set of a codimension one mass minimizer \( T_0 \) is covered by a countable collection of \( C^1 \) submanifolds of dimension \( m-7 \).
		\end{thm}
		
		A new and simpler proof of Simon's theorem has been recently found by Naber and Valtorta \cite{naber}. To the best of our knowledge, however, it is still an open question whether or not the remainder of the singularity set stratifies as lower-dimensional submanifolds.
		
		Surprisingly, codimension one mass minimizers do not have boundary singularities, as Hardt and Simon \cite{hardtsimon} established:

		\begin{thm}
			\label{thm:hardtsimon}
			If \( T_0\in \cal{I}_n(\R^{n+1}) \), \( \p T_0=[[S]] \) for some oriented embedded \( C^2 \) submanifold \( S\subset \R^n \), and the mass of \( T_0 \) is minimal among all integral currents with the same boundary, then there exists an open neighborhood \( V \) of \( S \) such that \( V\cap \supp(T_0) \) is an embedded \( C^{1,\alpha} \) hypersurface with boundary for all \( 0<\alpha<1 \).
		\end{thm}

		The story is more complicated and incomplete in higher codimension. Almgren in his 1700 page ``big regularity paper'' \cite{almgrenQ} proved the following theorem:

		\begin{thm}
			\label{thm:al}
			The singular set \( \cal{S}(T_0) \) of an \( m \)-dimensional mass-minimizing integral current \( T_0 \) in \( \R^n \) has Hausdorff dimension at most \( m-2 \).
		\end{thm}
		
		Again, this bound is sharp in codimension \( \geq 2 \) \cite{federer1965}. Chang \cite{chang} built upon Almgren's work to show that if \( m=2 \), then the singularity set consists of isolated branch points. More recently, in a series of papers \cite{delellisspadaro1,delellisspadaro2,delellisspadaro3,delellisspadaro4,delellisspadaro5}, De Lellis and Spadaro took on the monumental task of modernizing and simplifying Almgren's work. For an excellent overview of their approach, see \cite{delellisreview}.
		
\section{Non-Orientable Generalized Surfaces}
	Much of the above story can be repeated using chains with coefficients in a finite group, and in particular in \( \Z/2\Z \) to account for non-orientable surfaces. Fleming \cite{fleming} has a beautiful theory of flat chains with coefficients (see also \cite{ziemer},) the homology of which recovers the mod-\( p \) homology of the ambient space.
	
	However, soap films that occur in nature are not only non-orientable, but possess singularities such as triple junctions which are not amenable to mass-minimization. To ensure that the triple junction not be part of the algebraic boundary, one must assign one of the three surfaces a higher multiplicity. This in turn increases the total mass of the surface, and as a result triple junctions do not show up in solutions to the mass minimization problem.
	
	To get around this issue, there is a different approach one can take, and that is to ignore multiplicity when measuring area. Instead of minimizing mass, one can instead minimize \emph{\textbf{size}}, which for an integral current \( \sum_i k_i\int_{K_i} \) is the quantity \( \sum_i \H^k(K_i) \). The \( k \)-dimensional size of an arbitrary subset \( E\subset M \) is just \( \H^k(E) \). Note that the size of an integral current may be smaller than the size of its support. The primary difficulty with working with size is that unlike mass, it is not weakly lower semicontinuous. Extreme care must be taken with the minimizing sequence to account for this. The payoff is that size is better suited to the study of soap films than mass.
	
	\subsection{Reifenberg's 1960 Paper}
		\label{sub:reifenberg}
		The same year that Federer and Fleming's seminal paper appeared \cite{federerfleming}, Reifenberg published a work \cite{reifenberg} which dealt with the Plateau problem for non-orientable manifolds of arbitrary genus. This paper is also famous for a result that later became known as ``Reifenberg's disk theorem,'' which placed sufficient conditions on the approximate tangencies of a surface to guarantee that it was a topological disk. A set satisfying these conditions is now known as \emph{\textbf{Reifenberg flat.}} Reifenberg was well known for his prowess with tricky if not quirky estimates, and indeed his disk theorem did not disappoint: a condition involved in the statement required \( \e\leq 2^{-2000 n^2} \). Reifenberg's paper has sparked a number of subsequent results: Almgren provided a generalization to so-called ``elliptic integrands'' in \cite{almgrenannals}. Morrey generalized Reifenberg's result to ambient manifolds in \cite{morrey} (see also \cite{morreybook}.)
		
		Reifenberg's approach to proving his main theorem was built on work by Besicovitch and was purely set-theoretic, not involving any fancy machinery such as currents or varifolds. His main result was the following:
		
		Consider a finite collection \( A \) of pairwise disjoint Jordan curves in \( \R^3 \). A compact subset \( X \) of \( \R^3 \) is said to be a \emph{\textbf{surface spanning \( A \)}} if \( X \) can be written as an increasing union of manifolds \( X_i \) with boundary, such that for each \( i \), there exists a manifold \( Y_i \) with boundary \( A\cup \p X_i \) such that \( Y_i\to A \) in the Hausdorff distance. Reifenberg then proved:
		\begin{thm}
			There is a surface spanning \( A \) of least area.
		\end{thm}
		
		\begin{figure}
			\centering
			\subfloat[The boundary of a triple M\"{o}bius band]{%
				\includegraphics[width=.45\textwidth]{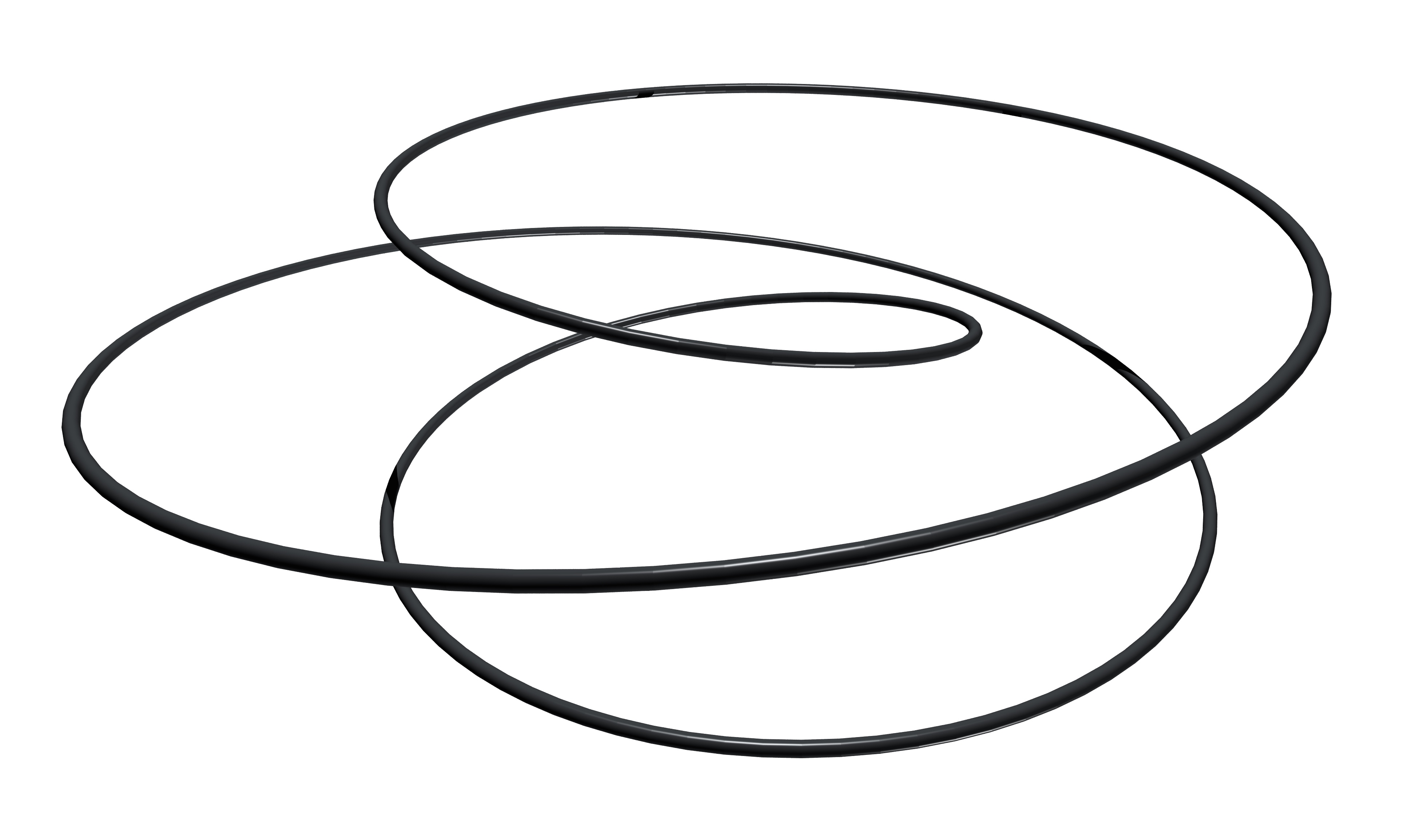}
				\label{fig:TripleWire}}
			\quad
			~
			\subfloat[A Triple M\"{o}bius band]{%
				\includegraphics[width=.45\textwidth]{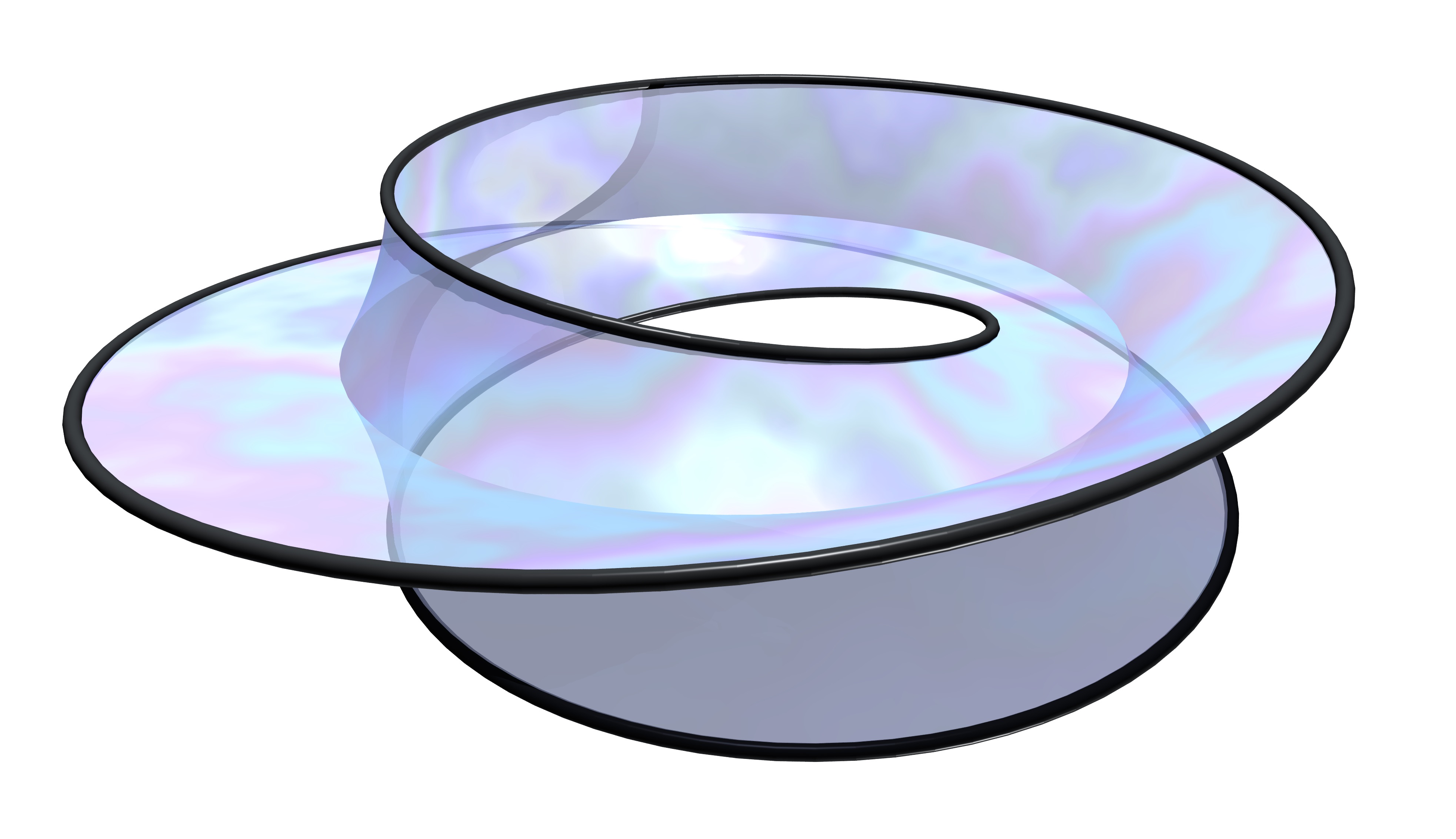}
				\label{fig:TripleMobiusStrip}}

			\subfloat[A non-orientable embedded surface]{%
				\includegraphics[width=.45\textwidth]{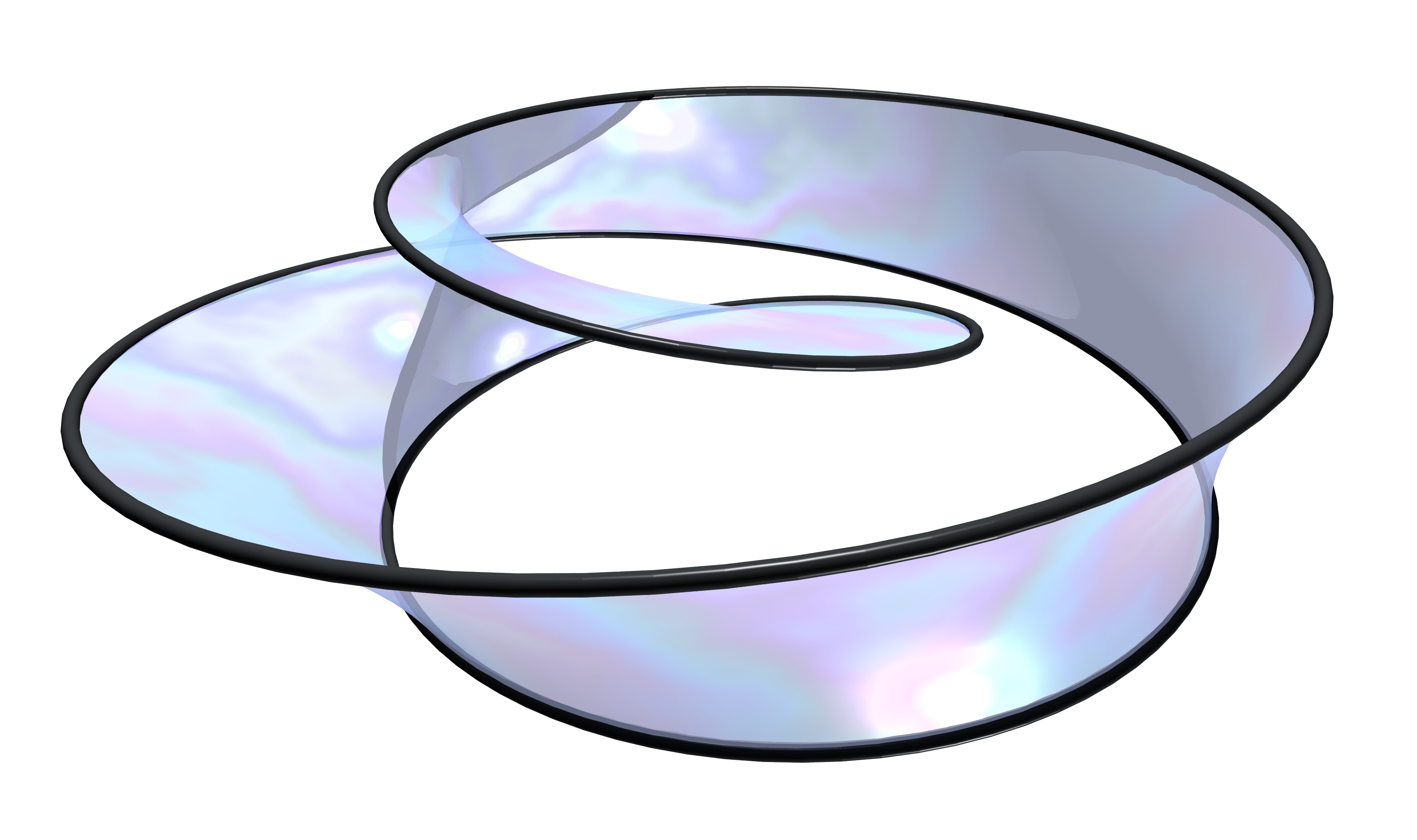}
				\label{fig:NonorientedEmbedded}}
			\quad
			~
			\subfloat[An immersed disk]{%
				\includegraphics[width=.45\textwidth]{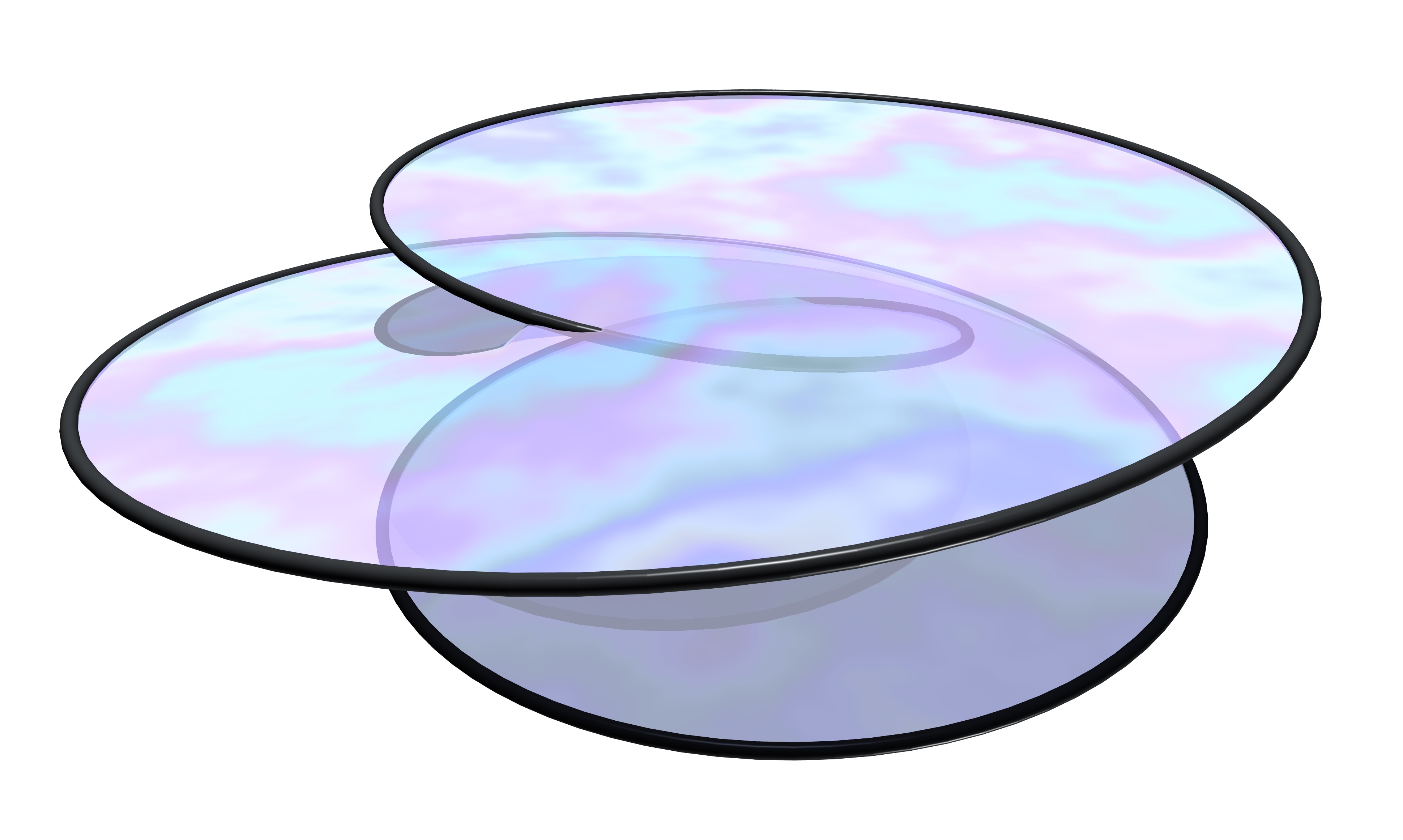}
				\label{fig:DiskImage}}
			\caption{\emph{Depending on the configuration of the boundary wire, either (b) or (c) can have smaller area.}}
			\label{fig:TripleMobius}
		\end{figure}
		
		His class of surfaces included those with infinite genus such as in Fleming's example \ref{fig:Fleming}, and non-orientable surfaces as well. However, it did not include soap-film type surfaces with triple junctions (Figure \ref{fig:TripleMobiusStrip}.) Reifenberg proved a secondary result which did minimize amongst this larger category, in the case that \( A\subset \R^n \) is homeomorphic to the \( (m-1) \)-sphere:

		\begin{thm}
			If \( A \) is a topological \( (m-1) \)-sphere in \( \R^n \), \( 2\leq m\leq n \) and \( \cal{G}^* \) is the collections of all compact sets \( X\supset A \) which do not retract onto \( A \), then there exists a set \( X\in \cal{G}^* \) with least \( m \)-dimensional Hausdorff spherical measure. Any such minimizer is locally Euclidian almost everywhere.
		\end{thm}
		
		Both of these theorems were special cases of a general result involving ``surfaces with algebraic boundary,'' which were defined by Reifenberg and developed by Adams in an appendix of \cite{reifenberg}.

		Let \( G \) be a compact\footnote{The exactness axiom is used in the proof and \v{C}ech homology only satisfies exactness when the coefficients are compact, so we will need this assumption} abelian group and suppose \( A \) is a compact subset of \( \R^n \) with \( \H^{m-1}(A)<\i \). Suppose \( L \) is a subgroup of the \( (m-1) \)-dimensional \v{C}ech homology \( \check{H}_{m-1}(A;G) \) of \( A \) with coefficients in \( G \). We say that a compact set \( X\supset A \) is a \emph{\textbf{surface with (algebraic) boundary \( \supset L \)}} if \( L \) is in the kernel of the inclusion homomorphism \( \iota_*: \check{H}_{m-1}(A;G)\to \check{H}_{m-1}(X;G) \). Reifenberg proved existence of a surface with algebraic boundary \( \supset L \) with least \( m \)-dimensional Hausdorff spherical measure. The case that \( G=\Z/2\Z \) implies Reifenberg's first theorem, and the case that \( G=S^1 \) implies, via a theorem of Hopf, the second.
		
\begin{figure}
			\centering	
			\includegraphics[scale=0.5]{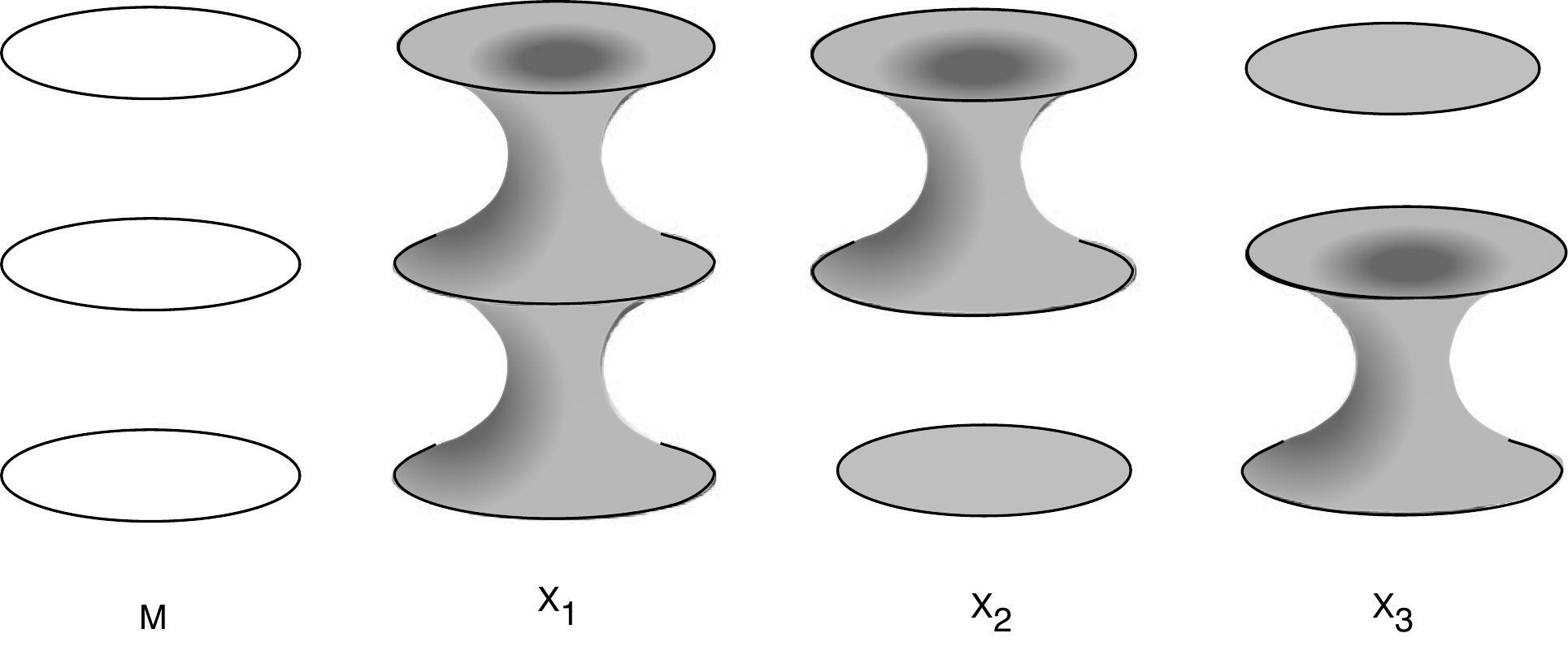}
			\caption{\emph{Three minimal surfaces spanning three circles}}
			\label{fig:ThreeCirclesProblem}
		\end{figure}

		A shortcoming of Reifenberg's theory is that for boundaries more general than a sphere, he did not defined a single, unifying collection surfaces with soap-film singularities. For example, consider the disjoint union of a disk and a circle in \( \R^3 \). There is no retraction to the pair of circles, yet we might not want to consider this as an admissible spanning set. As another example, consider the surfaces \( X_i \), \( i=1,2,3 \), in Figure \ref{fig:ThreeCirclesProblem}. Any one could be a surface with minimal area, depending on the distance between the circles, but a simple computation shows there is no non-trivial collection of Reifenberg surfaces which contains all three simultaneously. Thus, one would have to find an appropriate subgroup \( L \) which would produce the correct minimizer, and this task would change depending on the configuration of the circles in the ambient space.
		
		In a recent paper \cite{hpplateau,hpplateau2015}, the authors found a way around this problem using linking number to define spanning sets, and later in \cite{lipschitz}, using \v{C}ech cohomology in higher codimension. We also generalized Reifenberg's result so as to minimize a Lipschitz density functional.
		
		\subsubsection{Spanning sets via linking numbers}
			\label{ssub:spanning_sets_via_linking_numbers}
\begin{figure}
				\centering
				\includegraphics[width=.9\textwidth]{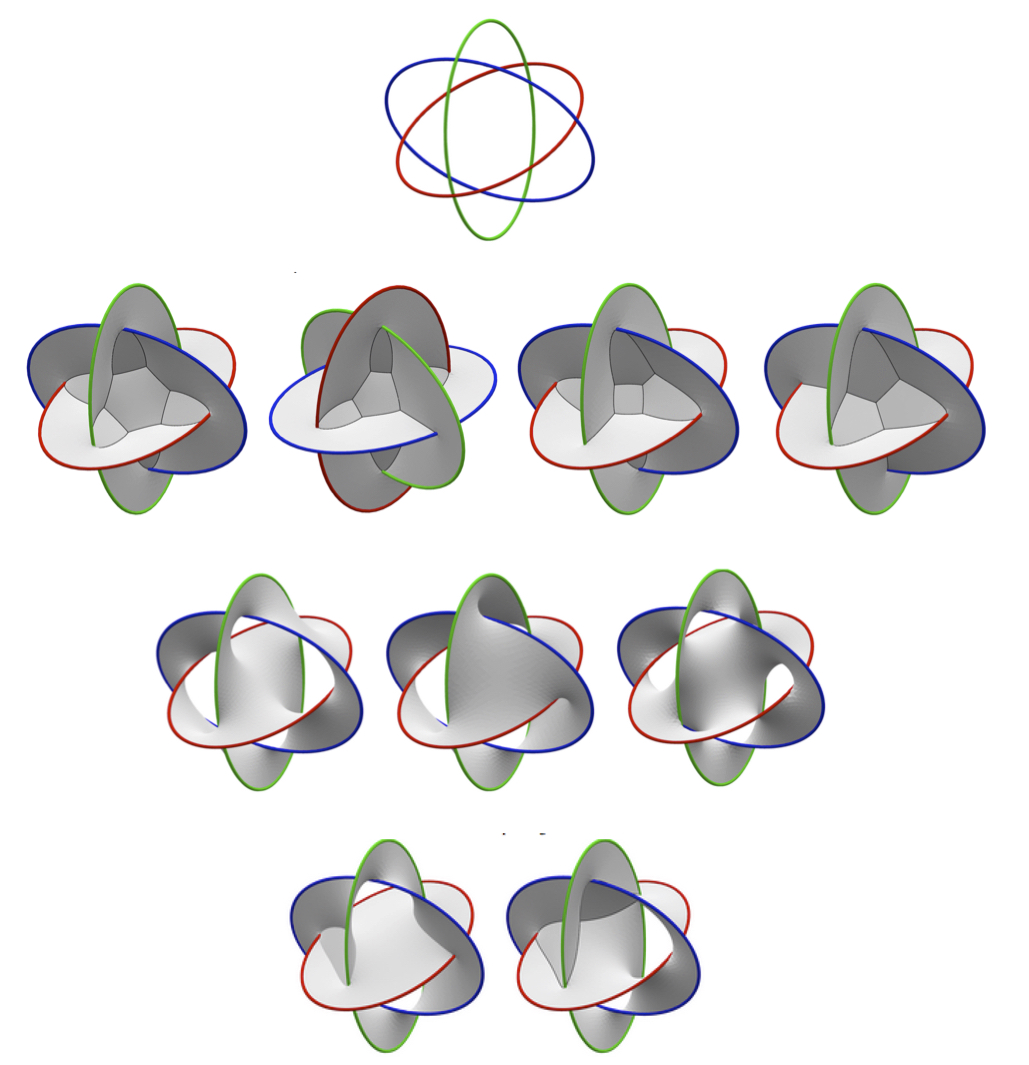}
				\caption{\emph{Each row of surfaces depicts a distinct type of minimal surface spanning the Borromean rings. A surface in the first row spans the Borromean rings using any linking test. That is, every simple link of any number of curves must meet the surface. In the second row, every simple link of one curve or all three curves must meet the surface. The third row has mixed types.}}
				\label{fig:BorromeanTests}
			\end{figure}

			\begin{defn}
				\label{def:linking}
				Suppose \( A \) is a \( (n-2) \)-dimensional compact orientable submanifold of \( \R^n \), \( n\geq 2 \). We say that a circle \( S \) embedded in \( \R^n\setminus A \) is a \emph{\textbf{simple link of \( A \)}} if the absolute value of the linking number \( L(S,A_i) \) of \( S \) with one of the connected components \( A_i \) of \( A \) is equal to one, and \( L(S,A_j)=0 \) for the other connected components \( A_j \) of \( A \), \( j\neq i \). We say that a compact subset \( X\subset \R^n \) \emph{\textbf{spans}} \( A \) if every simple link of \( A \) intersects \( X \) (See Figure \ref{fig:BorromeanTests}.)
			\end{defn}
			
			If \( A \) is a topological \( (n-2) \)-sphere then the set of spanning surfaces is the same as the collection \( \cal{G}^* \) above. Any orientable \( (n-1) \)-manifold with boundary \( A \) spans \( A \). The set \( A \) can be a frame such as the \( (n-2) \)-skeleton of an \( n \)-cube, in which case one can specify \( (n-2) \)-cycles which the simple links need to link. This procedure generalizes to higher dimension using linking spheres, or alternatively, via Alexander duality, to \v{C}ech cohomology. 

			This idea was first proposed for connected smooth boundaries in the existence paper of \cite{plateau10} which was followed by the more substantial \cite{hpplateau,hpplateau2015} which established lower semicontinuity of Hausdorff measure for a minimizing sequence \( X_k \to X_0 \) in codimension one and applied to any number of boundary components. One of us (HP) realized that linking number tests could be naturally viewed a cohomological spanning condition in higher codimension\footnote{Similar ``surfaces with coboundary'' were discovered independently by Fomenko \cite{fomenko}.}, and while we were writing this generalization \cite{lipschitz}, two papers appeared \cite{delellisandmaggi, ghiraldin} which built upon our linking number test for spanning sets.
			
			The cohomological spanning condition is stated as follows: if \( L \) is a \emph{subset} of the \( (m-1) \)-st (reduced) \v{C}ech cohomology group \( \check{H}^{m-1}(A;G) \) (\( G \) need not be compact,) we say that \( X\supset A \) is a \emph{\textbf{surface with (algebraic) coboundary \( \supset L \)}} if \( L \) is disjoint from the image of \( \iota^*: \check{H}^{m-1}(X;G)\to \check{H}^{m-1}(A;G) \).
			
			One of the primary benefits of using this definition over the covariant ``surface with algebraic boundary'' is that if \( A \) is an oriented manifold, then there is a natural choice for the subset \( L \), namely the collection \( L^{\Z} \) of those cocycles on \( A \) which evaluate to \( 1 \) on the fundamental cycle of a particular component of \( A \), and zero on the rest. By naturality of the Alexander duality isomorphism, the collection of surfaces with coboundary \( \supset L^{\Z} \) is equivalent in codimension one to the collection of compact sets which span \( A \) in the sense of linking number.

		Eight years after \cite{reifenberg} was published, Almgren proved an extension \cite{almgrenannals} of Reifenberg's theorem to prove the existence of surfaces which minimize not only area, but area weighted by a density function which is permitted to vary in both spacial and tangential directions, subject to an ellipticity condition. To discuss \cite{almgrenannals, almgren}, we will need to introduce rectifiable sets and varifolds, which we now define.

	\subsection{Rectifiable sets}
		\label{sub:rectifiable_sets_and_rectifiable_currents}
		A subset \( E \) of \( \R^n \) is \( m \)-\emph{\textbf{rectifiable}} if there exist a countable collection of Lipschitz maps \( \{f_i: \R^m\to \R^n\} \) such that the \( m \)-dimensional Hausdorff measure of \( E \setminus \cup_{i=0}^\i f_i(\R^m) \) is zero. If \( E \) is \( \H^m \) measurable and \( \H^m(E)<\i \), then the maps \( f_i \) can be taken to be \( C^1 \). Such sets are the higher dimensional analog of rectifiable curves. The defining property of a rectifiable set is that it is equipped with a unique ``approximate tangent \( m \)-plane'' almost everywhere, a consequence of Rademacher's theorem. If these approximate tangent spaces are equipped with an orientation, it becomes possible to integrate \( m \)-forms: An integer rectifiable current is just an integer weighted rectifiable \( E \) together with an orientation, i.e. an \( \H^m \) measurable field of \( m \)-vectors on \( E \) such that for almost every \( p\in E \), the \( m \)-vector at \( p \) is unit simple in the direction of the approximate tangent space at \( p \).

		A set \( F \) is \emph{\textbf{purely \( m \)-unrectifiable}} if \( \H^m(F\cap E)=0 \) for every \( m \)-rectifiable set \( E \). Every subset of \( \R^n \) with finite \( \H^m \) measure can be written, uniquely up to \( \H^m \) measure zero sets, as the disjoint union of a \( m \)-rectifiable set an a purely \( m \)-unrectifiable set. The beautiful Besicovitch-Federer structure theorem says that if \( F \) is purely \( m \)-unrectifiable, then for almost every \( m \)-plane \( V \) in the Grassmannian \( \Gr(m,n) \), the orthogonal projection of \( F \) onto \( V \) has \( \H^m \) measure zero.
		
		Morally, every subset of Euclidian space can be decomposed almost uniquely into a countable collection of \( C^1 \) submanifolds, and a remainder which casts no shadows.
		
	\subsection{Integral and Stationary Varifolds}
		
		Varifolds were first introduced by Young \cite{youngI, youngII} as ``generalized surfaces'' and developed by Young and Fleming \cite{flemingsurfaces, flemingyoung}. Fleming, who had been Young's student, in turn, taught Almgren what he knew \cite{flemingmemoir} when Almgren was a student at Brown 1958-62. Almgren took an interest in generalized surfaces and changed the name to ``varifolds,'' a mnemonic for manifolds in the calculus of variations. He produced a set of mimeographed notes \cite{almgrenmimeo} on varifolds that were circulated amongst his students but never published. Allard, who had been Fleming's student, produced the definitive reference on varifolds \cite{allard} in which he proved the compactness theorem for integral varifolds.

		\begin{defn}
			Let \( M \) be a smooth \( n \)-dimensional Riemannian manifold. \textbf{\emph{A \( k \)-varifold \( V \) in \( M \)}} is a Radon measure on the total space of the Grassmannian bundle \( \pi: \Gr\to M \), whose fiber above a point \( p\in M \) is the Grassmannian of un-oriented linear \( k \)-planes in \( T_p M \). The pushforward of \( V \) by \( \pi \) is denoted \( \|V\| \). The \emph{\textbf{mass of \( V \)}} is the quantity \( \|V\|(M) \). The \emph{\textbf{support of \( V \)}} is the support of the measure \( \|V\| \).
		\end{defn}

		For example, an embedded \( k \)-dimensional submanifold \( S\subset M \), together with a \( \H^k \)-measurable function \( \theta: S \to \R^+ \) determine a varifold \( V \) as follows: \[ V(A) := \int_{S \cap \{p: (p, T_p S) \in A\}} \theta(p) d\H^k(p). \]

		More generally, \( S \) can be replaced by a \( k \)-rectifiable set. In this case \( V \) is called a rectifiable \( k \)-varifold. If \( \theta \) takes integer values, the varifold is called an integral varifold. Integral varifolds are the non-orientable analogs of integral currents. There is no notion of integration of differential forms on a varifold, and unlike currents the space of varifolds does not possess a boundary operator. However, integral \( k \)-varifolds can be pushed forward by a Lipschitz map.
		
		Almgren saw these features as an advantage, for he wanted to model non-orientable surfaces and those with triple junctions. In \cite{varifolds} Almgren credits Federer-Fleming for proving the Plateau problem for mass minimization of oriented surfaces, and Reifenberg for size minimization of non-oriented surfaces ``subject to certain topological restraints.''  Almgren sought at this time to prove a mass minimization result for non-oriented surfaces. He did not specify what it meant for a varifold to span a given contour in \cite{almgrenmimeo} or \cite{varifolds}, but his focus at the time was on stationary varifolds where the definition seemed self-evident as we shall next see.

		The \emph{\textbf{first variation}} \( \d V \) of a compactly supported varifold \( V \) is a function which assigns to a smooth compactly supported vector field \( Y \) on \( M \) the rate of change of the mass of the pushforward of \( V \) by the time-\( t \) map of the flow of \( Y \) at \( t = 0 \). A varifold is \emph{\textbf{stationary}} if \( \d V = 0 \). 
		
		In \cite{almgrenmimeo}, Almgren proved the following theorem:
		\begin{thm}
			\label{thm:alm}
			Let \( M \) be a smooth compact \( n \)-dimensional Riemannian manifold. For each \( 0<k<n \) there exists a stationary integral \( k \)-varifold in \( M \).
		\end{thm}
		
		Allard in \cite{allard} proved a beautiful regularity result for such varifolds:
		\begin{thm}
			If \( V \) is a stationary integral \( k \)-varifold in a smooth compact \( n \)-dimensional riemannian manifold \( M \), \( 0<k<n \), then there is an open dense subset of the support of \( V \) which is a smooth \( k \)-dimensional minimal submanifold of \( M \).
		\end{thm}
 
		Even many years later, this theorem remains state-of-the-art in terms of what is known about the singularity set of stationary varifolds. For example, it is not known if the singularity set has zero Hausdorff measure in dimension \( k \). Indeed, regularity theory for stationary varifolds is still at an early stage, even compared to what is known about mass minimizing integral currents in higher codimension. It has been shown that the integral varifold
can be covered, up to a set of measure zero, by twice continuously
differentiable submanifolds of the same dimension, see \cite{schatzle}, \cite{menne}. See \cite{bombieri} who mentions this problem, as well as \cite{delellis} for a more detailed accounting of progress. 
    
	\subsection{Elliptic variational problems}
		\label{sub:subsection_Almgren_elliptic}

 %
 In \cite{almgrenannals} Almgren initiated the study of elliptic variational problems for non-orientable surfaces by providing the first definition of an elliptic integrand and a proof of regularity, depending on the degree of smoothness of the integrand.  His definitions and main regularity result follow: 
		
		Let \( A \) be a compact \( (m-1) \)-rectifiable subset of \( \R^n \) with \( \H^{m-1}(A)<\i \), \( G \) a finitely generated abelian group, and \( \sigma\in \check{H}_m(\R^n, A; G) \). We say a compact \( m \)-rectifiable set \( X\supset A \) is a \textbf{\emph{surface which spans \( \sigma \)}} if \( \sigma \) is in the kernel of the homomorphism on homology induced by the inclusion \( (\R^n,A)\hookrightarrow (\R^n,X) \).  
		
		A \emph{\textbf{\( C^k \) (resp. real analytic) integrand}} is a \( C^k \) (resp. real analytic) function \( f: \R^n \times \Gr(m,n) \to [a,b], \) where \( 0<a<b<\i \). We say \( f \) is \emph{\textbf{elliptic with respect to \( G \)}} if there exists a continuous function \( c:\R^n\to \R\cap \{t: t>0\} \) such that if \( D\subset \R^n \) is an \( m \)-disk, \( \tau\in \check{H}_m(\R^n, \p D; G)\setminus \{0\} \), and \( \tilde{D} \) is any surface which spans \( \tau \), then
		\[
		\int_{\tilde{D}} f(x, T_y \tilde{D}) \, d\H^m(y) - \int_D f(x, T_y D) \, d\H^m(y) \geq c(x) \left(\H^m(\tilde{D})-\H^m(D) \right)
		\]
		for all \( x\in \R^n \).
		
		\begin{thm}
			\label{theorem:almgren68}
			Let \( 3\leq k \leq \i \) and \( G\in \mathbf{G} \). If \( f \) is a \( C^k \) (resp. real analytic) integrand, elliptic with respect to \( G \), and  \( S \) is a surface that spans \( \sigma \) such that
			\[
			\int_S f(x, T_x S)\,d\H^m(x) \leq \int_T f(x, T_x T)\,d\H^m(x)
			\]
			for all surfaces \( T \) which span \( \sigma \), then  \( S \) is \( \H^m \) almost everywhere a \( C^{k-1} \) (resp. real analytic) submanifold of \( \R^n \).
		\end{thm}

   	   The authors proved there exists such a surface \( S \).  These results marked a significant advance\footnote{Readers should be warned that the existence portion of \cite{almgrenannals} contains a serious gap.  Briefly, a minimizing convergent sequence for a bounded elliptic integrand does not automatically yield a uniformly quasiminimal subsequence, but \cite{almgrenannals} assumes that it does. This is a critical part of the argument for existence of a minimizer (see \cite{elliptic} for a more detailed discussion.)   In \cite{almgren} spanning surfaces are chosen to be a priori uniformly quasiminimal so that the problem disappears. However, he was not able to prove a general existence theorem.} over Reifenberg's paper which only dealt with the functional \( f = 1 \), and bring Plateau's problem, which had grown well beyond the classical theory of minimal surfaces and the minimal surface equation, squarely back into the realm of PDE's. 

		\subsubsection{\( (f,\e,\d) \)-minimal sets}
   	    	\label{sub:_m_0_d_minimizing_sets}
   	    	In his memoir \cite{almgren} (see also \cite{almbulletin}), Almgren defined new classes of surfaces to model soap bubbles as well as many types of soap films.

			Fix \( A\subset \R^n \). If \( \phi:\R^n\to \R^n \) is Lipschitz, let \( W_\phi=\{x: \phi(X)\neq x \} \). If \( W_\phi\cup \phi(W_\phi) \) is disjoint from \( A \) and contained in a ball of radius \( \d \) for some \( 0<\d<\i \), we say that \( \phi \) is a \emph{\textbf{\( \d \)-deformation fixing \( A \)}}.

			Let \( 1\leq \g<\i \). A compact set \( X\subset \R^n \) with \( \H^m(X)<\i \) is \emph{\textbf{\( (\g,\d) \)-restricted with respect to \( A \)}} if \[ \H^m(X \cap W_\phi) < \g \H^m(\phi(X \cap W_\phi)) \] for all \( \d \)-deformations \( \phi \) fixing \( A \).
		
			If \( \e:[0,\i) \to [0,\i) \) with \( \e(0)=0 \) is a continuous non-decreasing function, the set \( X \) is called \emph{\textbf{\( (f,\e,\d) \)-minimal}} if in addition for every \( r \)-deformation \( \phi \) fixing \( A \), \( 0<r<\d \), \[ \int_X f(x,T_x X) \,d\H^m(x) \leq (1+\e(r))\int_{\phi(X)} f(x,T_x \phi(X))\, d\H^m(x). \]

			An important fact about \( (\g,\d) \)-restricted sets is that they are \( m \)-rectifiable (\cite{federer} 3.2.14(4)) and have both upper and lower bounds on density ratios. Almgren proved regularity results a.e for \( (f,\e,\d) \) minimal sets in  \cite{almgren}:
			\begin{thm}
				\label{thm:alreg}
				Suppose \( f \) is elliptic and \( C^3 \) and \( X \) is \( (f,\e,\d) \)-minimal with respect to \( A \) where \( \e \) satisfies \[ \int_0^1 t^{-(1+\alpha)} \e(t)^{1/2}dt < \i \] for some \( 0 \leq \alpha < 1 \). Then there exists an open set \( U \subset \R^n \) such that \( \H^m(X \setminus U) = 0 \) and \( X \cap U \) is a \( C^1 \) \( m \)-dimensional submanifold of \( \R^n \).
			\end{thm}
			Almgren states however that ``These hypotheses and conclusions, incidentally, do not imply that \( S \cap U \) locally can be represented as the graph of a function which satisfies any of the various Euler equations associated with \( f \). '' Indeed, any \( C^2 \) \( m \)-dimensional submanifold \( S \) with boundary is\footnote{Here and in the literature, \( M \) denotes the constant function \( 1 \).} \( (M,\e,\d) \) minimal with respect to \( \p S \) if \( \d \) is sufficiently small and \( \e \) is a linear map with large slope.
		
			For the three-dimensional case, Taylor \cite{taylor} relied upon Theorem \ref{thm:alreg} to prove a beautiful soap film regularity result for \( (M,\e,\d) \)-minimal sets. However, Morgan \cite{morgan} points out that the class of \( (M,0,\d) \)-minimal sets, taken over all \( \d > 0 \), is not compact. It remains an open problem of whether a smoothly embedded closed curve in \( \R^3 \) bounds a film with minimal area in the class of all \( (M,0,\d) \)-minimal sets.
			
			Another open problem motivated by \cite{delellisandmaggi, ghiraldin, hpplateau,hpplateau2015} is to prove the same regularity theorems as above in the case that \( \phi \) is also required to be uniformly close to a diffeomorphism.

\section{Variable boundaries}
	\label{sec:variable_boundaries}
	\subsection{Sliding boundaries}

		The notion of a sliding boundary has had a long history in the study of elasticity in mechanical engineering (see \S 24 of \cite{Podio-Guidugli}, for example.) David brought the attention of this problem to those in geometric measure theory \cite{david2014}. We shall mention a formulation of the problem found in \cite{delellisandmaggi, ghiraldin} which was influenced by \cite{david2014}. These works assume that the bounding set \( A \) has zero \( \H^m \) measure and often \( (m-1) \)-rectifiable, but others do not (see, e.g., \cite{coxjones, friedsequintubes}) as applications often require a large, and even rough, bounding set.
 
		\begin{defn}
			Let \( A \subset \R^n \) be compact and \( S_* \subset \R^n \setminus A \) be relatively compact. Let \( \S(A) \) denote the collection of Lipschitz maps \( \phi:\R^n \to \R^n \) such that there exists a continuous map \( \Phi:[0,1] \times \R^n \to \R^n \) with \( \Theta(1,\cdot) = \phi \), \( \Theta(0,\cdot) = Id \) and \( \Theta(t, A) \subset A \) for each \( t \in [0,1] \). Define \[ \cal{C}(A,S_*) = \{ S: S = \phi(S_*) \mbox{ for some } \phi \in \S(A) \} \] and call \( S_* \) a \emph{\textbf{sliding minimizer}} if \( \H^m(S_*) = \inf\{\H^m(S): S \in \cal{C}(A,S_*) \} \).
		\end{defn}
 
		Note that \( \cal{C}(A,S_*) \) does not form an equivalence class. It is not known if \( \cal{C}(A,S_*) \) is compact. However it can be shown with some assumptions on \( A \) (see \cite{delellisandmaggi} for codimension one, \cite{ghiraldin} for higher codimension) that if \( \{ S_k \}\subset\cal{C}(A,S_*) \) is a minimizing sequence, then the measures \( \H^m\lfloor_{S_k} \) converge weakly to a measure \( g \H^m\lfloor_{S_0} \) where \( S_0 \) is \( m \)-rectifiable and \( g \geq 1 \). In particular, \( \H^m(S_0) \leq \liminf \H^m(S_k) \). It is not known if \( S_0 \in \cal{C}(A,S_*) \), but \cite{delellisandmaggi} and \cite{ghiraldin} proved nonetheless that \( S_0 \) is a sliding minimizer.
		
		It is an open question if this result extends to Lipschitz or elliptic integrands. It is similarly open to prove the result if the aforementioned assumptions on \( A \) (e.g. \( \H^m(A)<\i \)) are removed.
		
	\subsection{Euler-Plateau problem}
		\label{sub:euler_plateau_problem}
		
		Mahadevan and Giomi \cite{mahadevan} proposed a type of Plateau problem in which a rigid boundary is replaced by a soft boundary such as a flexible wire. Specifically, in the language of Kirchhoff's theory of rods \cite{dill}, permissible boundaries are circular rods which resist bending yet are inextensible, unshearable, without intrinsic curvature, and without resistance to twisting about their centerlines. Mahadevan and Giomi formulated an energy functional which measured not only the area of a spanning surface, but also the energy of the boundary. The resulting Euler-Lagrange equations are equivalent, in the zero surface tension case, to those derived by Langer and Singer \cite{langersinger} (see \cite{friedchen}.) This minimization problem is called the \emph{\textbf{Euler-Plateau problem}} after Euler's study of column buckling, but might more appropriately be called the \emph{\textbf{Kirchhoff-Plateau problem}}.\footnote{The authors would like to thank Eliot Fried for helpful remarks.}
		
		Chen and Fried \cite{friedchen} rigorously derived the equilibrium conditions for the minimization problem, and provided geometric and physical interpretations of these conditions. Briefly, the surface on the interior must have zero mean curvature, and the boundary is required to bend elastically in response to a force exerted by the spanning film. The class of competitors for minimization are those surfaces which occur as images of the disk. However, since the boundary is permitted to vary, the maps cannot, in contrast to Douglas, Rad\'o and Courant, be assumed to be conformal. See also \cite{friedbiria, friedsequintubes}.

		These papers are closely related to earlier work by Bernatzki \cite{bernatzkielastic} and Bernatzki-Ye \cite{bernatzkiye}.

\section{Open problems}
	\label{sec:openproblems}
	
	Though we have labeled the following problems as ``open,'' some may have been solved without our knowledge. If you have solved one of these, please accept our apologies (and our congratulations!)
	
	\subsection{Classical minimal surfaces}
		The following problems are part of a longer list in \cite{meeksperez}. Let \( C \) be the space of connected, complete, embedded minimal surfaces and let \( P \subset C \) be the subspace of properly embedded surfaces.
		\begin{itemize}
			\item Isolated Singularity Conjecture (Lawson and Gulliver):  The closure of a properly embedded minimal surface in the punctured closed unit ball is a compact embedded minimal surface. 
			\item Convex Curve Conjecture (Meeks):  Two convex Jordan curves in parallel planes cannot bound a compact minimal surface of positive genus.
			\item \( 4 \pi \)-Conjecture  (Meeks, Yau, Nitsche): If \( \G \) is a simple closed curve in \( \R^3 \) with total curvature at most \( 4 \pi \), then \( \G \) bounds a unique compact, orientable, branched minimal surface and this unique minimal surface is an embedded disk. 
			\item Liouville Conjecture (Meeks): If \( M \in P \) and \( h:M \to \R \) is a positive harmonic function, then \( h \) is constant.
			\item Finite Genus Properness Conjecture (Meeks, P\'erez, Ros): If \( M \in C \) and \( M \) has finite genus, then \( M \in P \).
		\end{itemize}

	\subsection{Integral currents}
		These problems are adapted from a longer list in \cite{ambrosioreview}.
		\begin{itemize}
			\item Establish the uniqueness of tangent cones to an mass-minimizing current. Uniqueness for \( 2 \)-dimensional currents was proved in \cite{whitetangent}, and partial results in the general case in \cite{allardalmgren} and \cite{simona}.
			\item Does the singular set of a mass-minimizing current have locally finite \( \H^{m-2} \) measure? Chang \cite{chang} proved that it does if \( m = 2 \).
			\item Is the singular set of an mass minimizing current rectifiable? Does it have other geometric structure such as a stratification? See e.g. Theorem \ref{thm:sim95}.
		\end{itemize}

	\subsection{Reifenberg problems}
		Reifenberg posed ten open problems in \cite{reifenberg}. Three of particular interest are these:
		\begin{itemize}
			\item Let \( M \) be a manifold with boundary and \( D_k \) be discs with boundary \( \p M \). Let \( \mu \) be the infimum of the areas of discs with boundary \( \p M \). Suppose \( D_k \to M \) and \( \H^2(D_k) \to \H^2(M) = \mu \). Prove that \( M \) is a disc. Prove the same for \( m \)-dimensional disks.
			\item Generalize Theorem 2 of \cite{reifenberg} to the case where the boundary is any manifold. For example, let \( A \) be the \( 2 \)-torus and \( X \) the solid torus with a small interior ball removed. Then \( A \) is not a retract of \( X \), as one can deformation retract \( X \) onto the union of \( A \) and a transverse disk. If the torus is made to be narrower in some region, the solution to the generalized Theorem 2 in this case would be a transverse disk at the narrowest location.
			\item Find a class of surfaces which includes those such as the Adams example which retract onto their boundary, and also includes some class of deformations thereof; then prove a compactness theorem for such surfaces. Do those sets which do not admit a deformation retraction onto the boundary forms such a class?
		\end{itemize}

	\subsection{Elliptic integrands}
		\begin{itemize}
 			\item Show by example that interesting non-smooth solutions can arise which represent observed phenomena in nature if an elliptic integrand is not smooth.
			\item Prove a version of the main result in \cite{elliptic} for mass, instead of size, weighted by an elliptic density functional.
			\item What restrictions on the competing class of surfaces can be made that carry over to minimizing solutions? E.g., one can restrict the problem to graphs, disks, continuous embeddings, bordisms, topological type, etc. Each problem presents its own existence and regularity questions.
			\item  Axiomatic approach: Let \( \cal{S} \) be a collection of surfaces such that if \( S\in \cal{S} \) and \( \phi \) is a Lipschitz map fixing \( A \) which is \( C^0 \) close to a diffeomorphism, then \( \phi(S)\in \cal{S} \). What are minimal conditions needed on \( \cal{S} \) to guarantee existence of a minimizer in \( \cal{S} \) for an elliptic area functional?
			  See \cite{elliptic, delellisandmaggi, ghiraldin, DeLDeRGhi16}.
		\end{itemize}

	\subsection{Non-closed curves}
		\begin{itemize}
			\item Find models for surfaces spanning non-closed curves and prove a compactness theorem (Figure \ref{fig:Almgren}.) In \cite{lipschitz} we proposed using relative (co)homology as follows: If one replaces the boundary set \( A \) with a pair \( (A,B) \), the definition of a surface with coboundary can be repeated: A pair \( (X,Y)\supset (A,B) \) is a surface with coboundary \( \supset L \) if \( L \) is disjoint from the image of \( \iota^*:\tilde{H}^{m-1}(X,Y)\to \tilde{H}^{m-1}(A,B) \). To what extent can this be adapted if \( B \) is permitted to vary in some restricted fashion? See also \cite{drachmanwhite} and \cite{morgan} 11.3.
		\end{itemize}
	\subsection{Varifolds}
		\begin{itemize}
			\item Does a smoothly embedded closed curve in \( \R^3 \) bound a film with minimal area in the class of all \( (M,0,\d) \)-minimal sets?
			\item Does the singular set of a stationary varifold have measure zero?
		\end{itemize}

	\subsection{Dynamics and deformations}
		We close with four increasingly open-ended problems, ending with one contributed by the authors:
		\begin{itemize}

			\item Euler-Plateau for sliding boundaries: State and solve the Euler-Plateau problem for sliding boundaries. Not only is the bounding set \( B \) allowed to be flexible, but frontiers of solutions can slide around within \( B \) as it flexes.
			
			\item To what extent can mean curvature flow detect soap-film solutions including triple junctions and non-orientable surfaces, starting with a given spanning set? See \cite{brakke, whiteflow1,whiteflow}.

			\item The problem of lightning (Harrison and Pugh): Formulate a dynamic version of Plateau's problem which models the formation and evolution of branched solutions. Applications would be numerous: lightning, formation of capillaries, branches, fractures, etc. One should permit boundaries with higher Hausdorff dimension and solve the elliptic integrand problem, as solutions would be branched minimizers of the corresponding action principle, and thus should be highly relevant to physics.
		\end{itemize}

	\bibliography{bibliography.bib}{}
	\bibliographystyle{amsalpha}

\end{document}